\documentclass[a4paper, reqno]{amsart}

\usepackage[utf8]{inputenc}	
\usepackage[T1]{fontenc}
\usepackage{a4wide}
\usepackage{amsmath}		
\usepackage{enumerate}	
\usepackage{nicefrac}		
\usepackage{array}			
\usepackage{calc}				
\usepackage{flafter}	
\usepackage{esdiff}
\usepackage{pgffor}   
\usepackage{pgf}
\usepackage{mathabx}
\usepackage{multicol}
\usepackage{enumitem}

\usepackage
[
pdfauthor   = {Lisa Maria Kreusser},
pdftitle    = {Title},
pdfsubject  = {},
pdfkeywords = {},
bookmarks=true,
bookmarksopen=true,
pagebackref=false,
hyperindex=true,
colorlinks=true,
linkcolor=blue,
citecolor=blue,
filecolor=blue,
urlcolor=blue,
]
{hyperref}					
\usepackage{textcomp}		
\usepackage{color}		
\usepackage{graphicx}		

\usepackage{mathrsfs}		
\usepackage{graphicx}
\usepackage[caption=false,subrefformat=parens,labelformat=parens]{subfig}

\usepackage{bbm}

\newcommand\N{\mathbb{N}}
\newcommand\R{\mathbb{R}}

\newcommand\D{\mathcal{D}}
\newcommand\E{\mathcal{E}}

\newcommand\bl{\left(}
\newcommand\br{\right)}

\newcommand*\di{\mathop{}\!\mathrm{d}}

\newcommand\bv{\operatorname{BV}}
\newcommand\tv{\operatorname{TV}}

\renewcommand\epsilon{\varepsilon}

\renewcommand\theta{\vartheta}

\newtheoremstyle{mytheoremstyle} 
{6pt}                  
{6pt}                    
{\itshape}               
{}					
{\bf}                  
{}                    
{.5em}                 
{} 
\theoremstyle{mytheoremstyle}
\newtheorem{satz}{Satz}[section]
\newtheorem{lemma}[satz]{Lemma}
\newtheorem{proposition}[satz]{Proposition}
\newtheorem{corollary}[satz]{Corollary}
\newtheorem{theorem}[satz]{Theorem}

\theoremstyle{remark}
\newtheorem{ass}[satz]{Assumption}

\graphicspath{{./../}}

\numberwithin{equation}{section}
\usepackage{cleveref}

\makeatletter
\def\namedlabel#1#2{\begingroup
	#2%
	\def\@currentlabel{#2}%
	\phantomsection\label{#1}\endgroup
}
\makeatother

\title{On anisotropic diffusion equations for label propagation}
\author{Lisa Maria Kreusser}
\address{Lisa Maria Kreusser, Department for Applied Mathematics and Theoretical Physics, University of Cambridge, Wilberforce Road, Cambridge CB3 0WA, UK}
\author{Marie-Therese Wolfram}
\address{Marie-Therese Wolfram, Mathematics Institute, University of Warwick, CV47AL Coventry, UK and RICAM, Austrian Academy of Sciences, Altenbergerstr. 69, 4040 Linz, Austria}

\begin{document}
	
	\begin{abstract}
	In many problems in data classification one wishes to assign labels to points in a point cloud with a certain number of them being already correctly labeled. In this paper, we propose a microscopic ODE approach, in which information about correct labels is propagated to neighboring points. Its dynamics are based on alignment mechanisms, which are commonly used in large interacting agent systems in consensus formation. We derive the respective continuum description, which corresponds to an anisotropic diffusion equation with reaction term. Solutions of the continuum model on the bounded domain inherit certain properties of the underlying point cloud. We discuss these analytic properties and exemplify the results with micro- and macroscopic simulations.
	\end{abstract}
	
	\maketitle

	\section{Introduction}

	Large interacting agent systems have been used successfully to model consensus formation and collective dynamics in the life and social sciences, such as 
	fish schools, bird flocks or large societies. 
In these systems imitation, alignment and other interactions lead to the formation of 
complex patterns and stationary states such as flocks or clusters. This collective and self-organized behavior initiated a lot of research
in the applied maths community, but also inspired the development of 
many algorithms in data science, machine learning and robotics.  \\
In this paper we will apply ideas from consensus and collective dynamics in large interacting agent systems for labeling and classification problems. In particular we wish 
to assign labels to a given set of $n$ data points,
some of them being already correctly matched to labels - a situation known as semi-supervised learning. We assume that
information about correct labels is then propagated to unlabeled points using ideas from large interacting agents systems, see
\cite{motsch2014}.
Hereby interactions depend on the closeness of points - as in interacting agents systems, particles only interact with others within a certain radius $\epsilon > 0$. 
A generic form of such models is as follows:
let $u_i = u_i(t)$, $i=1, \ldots, n$, denote a characteristic of agent $i$, for example its opinion, and $w_{ij} = \eta(x_i, x_j, \epsilon)$ the influence that agent $i$ has 
on $j$. Then agents interact via
\begin{align}\label{eq:consensus}
\frac{\di u_i}{\di t} = \sum_{j \neq i} w_{ij} (u_j - u_i)  \text{ for } i = 1, \ldots, n.
\end{align}
Several well-known models fall into this category - for example  Krause's opinion formation model, 
see \cite{K2000}, or consensus formation models in large interacting agent systems, see \cite{motsch2014}. But also 
the well-known Cucker-Smale model, see \cite{CS2007}, can be written as
\eqref{eq:consensus}. Hereby $u_i$ corresponds to the particle's velocity, with weights $w_{ij}$ depending on the position.
Related interacting agent systems
with an anisotropic interaction or influence function have only recently been considered, for example in the context of fingerprint formation, 
see \cite{DGHKS2017}. \\
We  will consider \eqref{eq:consensus} with an additional regularization term, which enforces binary labels of points through a double well potential. 
Furthermore we assume that already labeled points only influence others via interaction  - a situation similar to interactions with leaders 
in opinion formation, see \cite{DMPW2009}.\\
Models of type \eqref{eq:consensus} have been used not only in biology, but are also well known in classification and graph-partitioning problems. 
Hereby agents relate to data points or vertices of a weighted graph; the interaction function $w_{ij}$ to the edge weights. System \eqref{eq:consensus} then
relates to the graph Laplacian. Its continuum limit has been studied in different settings - going back to the works of Hein et al. \cite{HAL2007}, or more recent
improved results using methods such as $\Gamma$-convergence, see \cite{C2018,CST2020, CT2020, TS2016, ST2019}.
A detailed analysis of the eigenmaps and eigenvectors of the corresponding weighted graph Laplacians in the case of 
two nearly separated clusters,
was published quite recently by Hoffmann et al, see \cite{HHOS2019}. Hereby they use the Fiedler vector - the eigenvector of the second smallest 
eigenvalue of the graph Laplacian - to partition the graph. The identified continuum operators depend on the particular weighting of the graph
Laplacian - they correspond to the following  class of anisotropic elliptic operators 
\begin{align}\label{eq:adiff}
\mathcal{L}(u) = \frac{1}{\rho^p} \nabla \cdot \left(\rho^q \nabla \left(\frac{u}{\rho^r}\right)\right),
\end{align}
where $\rho$ denotes the distribution of points with parameters $p,~q,~r \in \mathbb{R}$ fixed. \\
\noindent Graph Laplacians in connection with suitable regularization terms have been investigated in various applications in data science and imaging in the last years.
The idea and notion of label propagation was first explored in the machine learning community, see \cite{ZG2002, ZGL2003}. Subsequent works include for
example van Gennip et al.\ who studied the asymptotic behavior of the graph Laplacian and a double well potential. This approach corresponds to the graph-based formulation
of the well-known Allen-Cahn equation, see \cite{gennip2012}. Merkurjev et al.\ \cite{MKB2013} used related ideas, presented in \cite{BF2012}, to propose computationally 
efficient algorithms to calculate eigenvalues and eigenvectors of graph Laplacians 
and improve the performance of existing methods, such as the Merriman-Bence-Osher algorithm. Note that Garcia-Cordona et al.\ proposed a related microscopic approach in \cite{GMBFP2014}. Hereby they enforce known correct labels by an additional term in their energy functional. However, our approach allows us to derive respective continuum description and analyse the structure of solutions.\\

We will present two alternative formal derivations of the limiting continuum model for the proposed microscopic labeling model. The first approach investigates the
behavior of the respective energy functionals,
which accounts for already correctly labeled points and their specific interactions in the simultaneous limit $n\rightarrow \infty$ and $\epsilon\rightarrow 0$.
In the second approach we consider the limiting problem starting from the proposed ODE system instead of the energy functional and hence no analysis of the spectrum is required. 
The proposed arguments allow us to deduce the respective continuum equation for 
semi-supervised learning,  in which correctly labelled points correspond to Dirichlet boundary conditions. In both approaches the derived interaction operators are consistent 
with \eqref{eq:adiff}, providing therefore an alternative viewpoint and derivation. \\
The contributions of this paper can be summarized as follows:
\begin{enumerate}
	\item \textit{Model development:} We propose a different perspective on label propagation, which connects consensus dynamics in large interacting systems with 
	graph-Laplacians in data science. The respective methodologies are well established in their fields, the proposed model adapts them in the context of semi-supervised 
	learning and establishes
	new connections.
	\item \textit{Continuum models:} We present two alternative formale derivations of the respective continuum model using an energy based approach and a discrete to continuum limit. 
	\item \textit{Quantitative behavior:} We present several analytic results to characterize the behavior of solutions to the continuum problem with suitable boundary conditions on bounded
	domains. 
	These results show that
	solutions to the continuum problem inherit the structure of the underlying data cloud.
	\item \textit{Computational experiments:} We illustrate and exemplify the structure of solutions to the micro- as well as macroscopic equation with various computational experiments.\\
\end{enumerate}

This paper is organized as follows. We introduce the microscopic model for label propagation in Section \ref{s:micmod}, before continuing with the derivation of the 
macroscopic description from the energy formulation and the formal derivation of the discrete to continuum limit in Section \ref{s:mictomac}.  We discuss the properties of solutions to the 
continuum model in Section \ref{s:pde} and conclude with computational experiments in Section \ref{s:numerics}.

\section{Description of the discrete model}\label{s:micmod}
In the following, we propose a model for the propagation of labels of a given point cloud.
We start by introducing the proposed label propagation model on the microscopic level. We consider 
a point cloud $$V=\{X_1,\ldots,X_{n+m}\}$$ of data points $X_i$ in some bounded domain $\Omega\subset \R^d$ for $d\geq 1$ and $m,n\in \N$.  
In the following, we assume that the first $n$ points $X_i$ for $i=1,\ldots, n$, are not yet labeled,  
while the points $X_i$ for $i= n+1,\ldots, n+m$ are equipped with labels $L=\{l_1,\ldots,l_k\}$ for some $k\ll n$ with $l_1<\ldots<l_k$. 
We denote the set of points associated to each label by $V_i'$ and assume that $V_1',\ldots,V_k'\subset V$ such that $V_i'\cap V_j'=\emptyset$ for $i,j\in\{1,\ldots,k\}$ with $i\neq j$.\\
We propagate information from close correctly labeled points to neighboring unlabeled ones. 
This assumption is well established in collective dynamic models, in which interactions depend on the metric or topological distance, see for example \cite{CYP2017, H2013}. We relate the
closeness of two data points to their metric distance and from now on consider the topology of the underlying respective graph $G$ and its weights $w$. 
This leads to a weighted graph $G=(V,E,w)$, where $V$ corresponds to the set of data points or agents, $E$ the connecting edges and $w$ the fixed interaction rates or weights.
We assume that $G$ is connected. This implies that each vertex is connected to at least another one and hence for 
each $i\in V$ there exists $j\in V\backslash \{i\}$ such that $w_{ij}\neq 0$. In the following, we also assume that the weights are symmetric, i.e.\ $w_{ij}=w_{ji}$, and 
we consider two clusters, i.e. $k=2$.

Let $\gamma>0$ and $\kappa>0$ be given constant and let $W\colon \R\to \R$ denote a double well potential with wells at 
$l_1=-1$ and $l_2=1$, for example\ $W(x)=(x^2-1)^2$. Given a graph $G$ with weights $w_{ij}$ and $V'=V_1'\cup V_2'$, we assume that labels $u_i$ for $i\in V$ of unlabeled
points change
according to
\begin{align}\label{eq:modelorig}
\begin{split}
\frac{\di u_i}{\di t} &=\gamma \sum_{j\neq i} w_{ij} (u_j -u_i)-\kappa W'(u_i), \quad i=1,\ldots,n.
\end{split}
\end{align} 
We can extend the sum in \eqref{eq:modelorig} to the entire set of points $V$. Therefore
\begin{align*}
\begin{split}
\frac{\di u_i}{\di t} &=\gamma \sum_{j=1}^{n+m} w_{ij} (u_j -u_i)-\kappa W'(u_i), \quad i=1,\ldots,n.
\end{split}
\end{align*} 
Since we are interested in the large number particle limit, we rescale the amplitude of the interaction term as follows 
\begin{subequations}
	\label{eq:micro}
	\begin{align}\label{eq:model}
	\begin{split}
	\frac{\di u_i}{\di t} &=\frac{\gamma}{ n+m} \sum_{j=1}^{ n+m} w_{ij} (u_j -u_i)-\kappa W'(u_i), \quad i=1,\ldots,n.
	\end{split}
	\end{align} 
	Labels of points in $V_1'$ and $V_2'$ do not change, but influence the dynamics via the interaction term in \eqref{eq:model}. Hence the following equalities have to be added to \eqref{eq:model} 
	\begin{align}\label{eq:bc}
	\begin{split}
	u_i=-1, \quad i\in V_1',\qquad
	u_i=1, \quad i\in V_2'.
	\end{split}
	\end{align}
	Furthermore we supplement the system with the initial label distribution 
	\begin{align}\label{eq:initialdiscrete}
	u_i(0)=-1, \quad i\in V_1', \qquad 
	u_i(0)=1, \quad i\in V_2',\qquad
	u_i(0)=0, \quad i\in V\backslash V'. 
	\end{align}
\end{subequations}
Then the corresponding energy functional  is given by
\begin{align}
\begin{split}
\label{eq:energyorig}
E(u_1,\ldots,u_n)&=\frac{\gamma}{n(n+m)}  \sum_{i=1}^n\sum_{j=1}^{n} w_{ij}(u_j-u_i) \\
&\quad +\frac{\gamma}{n(n+m)}  \sum_{i=1}^n\sum_{j=n+1}^{n+m} w_{ij}(u_j-u_i)+\frac{\kappa}{n} \sum_{i=1}^n W(u_i),
\end{split}
\end{align}
where we rescaled by $n$ due to the summation over $i=1,\ldots,n$. The slightly unusual form of the energy $E$ stems from dependence of $E$ on $u_i$ for $i=1,\ldots,n$ 
while $u_i$ for $i=n+1,\ldots,n+m$ are fixed by conditions \eqref{eq:bc}.

\section{Derivation of the macroscopic model}\label{s:mictomac}
We continue by presenting two different approaches to derive the continuum description 
in the limit $n\rightarrow \infty$ and $\epsilon \rightarrow 0$ - the first one is a formal argument based on the 
energy formulation, while the second one is derived  from the respective ODE system.

Let $\rho \in \mathcal{P}(\Omega)$ denote the density of independent and identically distributed  points $V=\{X_1,\ldots, X_{n+m}\}$
and $\mathcal P(\Omega)$ corresponds to the space of probability measures on $\Omega$. 
The points $X_1,\ldots, X_{n+m}\in \R^d$ define the associated empirical measure $\rho_{n}(x)=\frac{1}{n+m}\sum_{i=1}^{n+m}\delta_{X_i} (x)$. 
Furthermore we introduce the function $u\colon \Omega \to \R$ satisfying $u(X_i)=u_i$ for all $i =1,\ldots n+m$. 
To rescale the discrete weights we consider the function $\bar \eta\colon \R^d\to [0,\infty)$. 
More precisely, we assume that the kernel $\bar \eta$ is isotropic and given by the radial profile $ \eta \colon [0, \infty) \to [0, \infty)$, 
i.e.\ $\bar \eta(x) = \eta(|x|)$, satisfying
\begin{enumerate}
	\item  $ \eta(0) > 0$ and $\eta$ is continuous at 0. 
	\item $ \eta$ is non-increasing.
	\item $\int_0^\infty  \eta(s ) s^{d+1} \di s<+\infty$. 
\end{enumerate}
Due to Assumption (3) and the radial symmetry of $\eta$ we define
\begin{align}\label{eq:sigmaeta}
\sigma_\eta:=\frac{1}{2}\int_{\R^d} \eta(|x| ) |x\cdot r|^{2} \di x<+\infty,
\end{align}
for any normal vector $r\in \R^d$. 
Note that the assumptions on $\eta$ are not restrictive and include a broad class of  kernels  such as Gaussian kernels and discontinuous kernels 
like $\eta(s)=\mathbbm{1}_{[0,R]}(s)$, where $\mathbbm{1}_{[0,R]}$ denotes the indicator function on $[0,R]$ for some $R>0$. 
We will later consider the appropriately rescaled function $\eta_\epsilon(s)=\epsilon^{-d}\eta(\frac{s}{\epsilon})$ in the limit. We assume that significant weight is given to edges connecting points up to distance $\epsilon$ and that the weights of the graph are given by $w_{ij}=\eta_\epsilon(|X_i-X_j|)$. With the limit $n\to \infty$ in mind, this suggests to consider a given scaling $\epsilon$ with respect to $n$. In addition, the interaction term in \eqref{eq:model} has to be rescaled by $\epsilon^2$, resulting in
\begin{align}\label{eq:modelweight}
\begin{split}
\frac{\di u_i}{\di t} &=\frac{\gamma}{\epsilon^2  (n+m)} \sum_{j=1}^{ n+m} \eta_\epsilon(|X_i-X_j|) (u_j -u_i)-\kappa W'(u_i), \quad i=1,\ldots,n,
\end{split}
\end{align} 
while  condition \eqref{eq:bc} implies that
\begin{align}\label{eq:modellabels}
\frac{\di u_i}{\di t}=0,\quad i= n+1,\ldots,n+m.
\end{align}
Similarly, we rescale the interaction term in the energy functional \eqref{eq:energyorig} by $\epsilon^2$ which yields
\begin{align}\label{eq:energy}
\begin{split}
E_{\epsilon,n}&(u_1,\ldots,u_n)=\frac{\gamma}{4\epsilon^2 n(n+m)}  \sum_{i=1}^n\sum_{j=1}^{n} \eta_\epsilon(|X_i-X_j|)(u_j-u_i)^2\\&\quad +\frac{\gamma}{2\epsilon^2 n(n+m)}  
\sum_{i=1}^n\sum_{j=n+1}^{n+m} \eta_\epsilon(|X_i-X_j|)(u_j-u_i)^2+\frac{\kappa}{n} \sum_{i=1}^n W(u_i).
\end{split}
\end{align}
Note that \eqref{eq:modelweight} corresponds to the $L^2$ gradient flow of the energy \eqref{eq:energy}.

\subsection{Properties of the discrete system for $\epsilon>0$}	
We continue by discussing the existence of solutions of the discrete system for $\epsilon>0$ under the assumption that the initial condition satisfies
\begin{align}\label{eq:energyinitbound}
\sup_{\epsilon>0}\sup_{n\in \N} E_{\epsilon,n} (u_1(0),\ldots,u_n(0))<\infty.
\end{align}
Note that we will use the shorthand notation $u(t)=(u_1(t),\ldots,u_n(t))$ in the following.
\begin{proposition}\label{prop:discretesol}	
	Let $\epsilon>0$ and $n\in\N$ be given and suppose that the initial condition satisfies \eqref{eq:energyinitbound}. Then, the regularized discrete system \eqref{eq:modelweight} with  condition \eqref{eq:bc} and initial data \eqref{eq:initialdiscrete} has a unique solution on $[0,T]$ which is uniformly bounded with respect to $t\in [0,T]$, $n\in\N$ and $\epsilon>0$. In particular, 
	\begin{align}\label{eq:energyeqdiscrete}
	E_{\epsilon,n}(u_1(t),\ldots,u_n(t))+\frac{1}{n} \sum_{i=1}^n \int_0^t \bl\frac{\di u_i(s)}{\di t}\br^2 \di s=E_{\epsilon,n}(u_1(0),\ldots,u_n(0)).
	\end{align}
\end{proposition}
\begin{proof}
	Since $W'(u_i)$ is monotonically increasing in $u_i$ outside the compact interval $[-1,1]$, Picard Lindel\"of ensures the existence of a unique, 
	continuously differentiable solutions to \eqref{eq:modelweight}.
	The boundedness of the solution follows from an energy identity. Using \eqref{eq:modelweight} and \eqref{eq:modellabels}, we obtain
	\begin{align*}
	\frac{\di }{\di t}E_{\epsilon,n}(u_1,\ldots,u_n)= -\frac{1}{n} \sum_{i=1}^n \bl\frac{\di u_i}{\di t}\br^2.
	\end{align*}
	This yields \eqref{eq:energyeqdiscrete}.
	The uniform boundedness of $E_{\epsilon,n} (u(t))$ with respect to $\epsilon$ and $n$ follows from the uniform boundedness of $E_{\epsilon,n} (u(0))$.
\end{proof}	

\subsection{Energy approach}\label{sec:energyapproach}
First we consider the formal limit of the energy functional $E_{\epsilon,n}$ as $n\to \infty$ and $\epsilon\to 0$,  before presenting an ODE-based limit in 
the next section. \\
Let \eqref{eq:initialdiscrete} satisfy \eqref{eq:energyinitbound}. Energy \eqref{eq:energy} can be rewritten as
\begin{align*}
E_{\epsilon,n}(u)
&=\frac{ \gamma}{4 \epsilon^2 n(n+m)}\sum_{i =1}^{n+m} \sum_{j =1}^{n+m} \eta_\epsilon(|X_i-X_j|)(u_j-u_i)^2\\&\quad-\frac{ \gamma}{4 \epsilon^2 n(n+m)}\sum_{i=n+1}^{n+m}\sum_{j =n+1}^{n+m} \eta_\epsilon(|X_i-X_j|)(u_j-u_i)^2   +\frac{\kappa}{n} \sum_{i=1}^{n+m} W(u_i)
\end{align*}
where $W(u_i)=0$ for $i=n+1,\ldots,n+m$ by the definition of $W$.
Since $\frac{1}{n}=\frac{1}{n+m}+\frac{m}{n(n+m)}$,
we obtain 
\begin{align*}
E_{\epsilon,n}(u)
&=\frac{ \gamma}{4 \epsilon^2 (n+m)^2}\sum_{i =1}^{n+m} \sum_{j =1}^{n+m} \eta_\epsilon(|X_i-X_j|)(u_j-u_i)^2+\frac{\kappa}{n+m} \sum_{i=1}^{n+m} W(u_i)\\&\quad +\frac{\kappa m}{n(n+m)} \sum_{i=1}^{n+m} W(u_i)+\frac{\gamma m}{4 \epsilon^2 n(n+m)^2}\sum_{i =1}^{n+m} \sum_{j =1}^{n+m} \eta_\epsilon(|X_i-X_j|)(u_j-u_i)^2\\&\quad-\frac{ \gamma}{4 \epsilon^2 n(n+m)}\sum_{i=n+1}^{n+m}\sum_{j =n+1}^{n+m} \eta_\epsilon(|X_i-X_j|)(u_j-u_i)^2.
\end{align*}
Note that all terms of $E_{\epsilon,n}$ except the first and second term vanish  in the limit $n\to \infty$ due to the uniform boundedness of $u_i$. 
By considering the empirical measure $\rho_n$ and
assuming that $\rho$ is bounded in $L^\infty(\Omega)$ as well as $\rho_n\to \rho$ weakly as $n\to \infty$, we obtain
\begin{align}\label{eq:energyresc}
\begin{aligned}[t]
E_\epsilon(u)=\lim_{n\to\infty} E_{\epsilon,n}(u)&=\frac{ \gamma}{4}\int_{\Omega}\int_{\Omega}\epsilon^{-2} \eta_\epsilon(|y-x|)\bl u(y)-u(x)\br^2\rho(x)\rho(y)\di y \di x\\
&\quad +\kappa \int_{\Omega}W(u(x))\rho(x)\di x.
\end{aligned}
\end{align}
Since $\Omega$ is a bounded domain, the set $\frac{\Omega-x}{\epsilon}$ for $x\in \Omega$ contains a neighborhood of the origin whose diameter increases as $\epsilon>0$ decreases. 
We recall that the radially symmetric mollifier $\eta$ has a compact support, given by the ball $B_r(0)$ for some $r>0$.   In particular, for $x\in \Omega$ given, there exists $\epsilon_0(x)>0$ such that $B_r(0)\subset \frac{\Omega-x}{\epsilon}$ for any $\epsilon\in(0,\epsilon_0(x))$.
Using a change of variable, we can consider the limit $\epsilon\to 0$ of the first term in \eqref{eq:energyresc}:
\begin{align*}
&\lim_{\epsilon\to 0}\frac{ \gamma}{4}\int_{\Omega}\int_{\Omega}\epsilon^{-2} \eta_\epsilon(|y-x|)\bl u(y)-u(x)\br^2\rho(x)\rho(y)\di y \di x \\
&=\lim_{\epsilon\to 0}\frac{\gamma}{4} 	\int_\Omega \int_{\frac{\Omega-x}{\epsilon}\cap B_r(0)} \eta(|z|) \bl \frac{u(x+\epsilon z) -u(x)}{\epsilon}\br^2 \rho(x)\rho(x+\epsilon z)\di z \di x\\
&=\frac{\gamma}{4} 	\int_\Omega \int_{ B_r(0)} \eta(|z|) \bl \nabla u(x) \cdot z\br^2 \rho^2(x)\di z \di  x\\
&=\frac{\gamma}{4} \sum_{i=1}^d	\int_\Omega \int_{ B_r(0)} \eta(|z|) \bl \partial_{x_i} u(x)\br^2 z_i^2 \rho^2(x)\di z \di  x=\frac{\gamma \sigma_\eta}{2}  \int_\Omega  \rho^2(x) |\nabla u(x)|^2 \di x,
\end{align*} 
where we used that the integral of an odd function on a symmetric domain vanishes in the second to the last equality 
and definition of $\sigma_\eta$ \eqref{eq:sigmaeta} in the last.
Then the resulting limiting energy is given by
\begin{align}\label{eq:energycont}
\E(u)=\frac{\gamma  \sigma_\eta}{2}\int_{\Omega} \rho^2 |\nabla u|^2\di x +\kappa \int_{\Omega}\rho W(u)\di x.
\end{align}
Its associated $L^2$-gradient flow is given by
\begin{align*}
\partial_t  u=\gamma  \sigma_\eta \nabla\cdot (\rho^2 \nabla u )-\kappa \rho W'(u),\quad x\in \Omega\backslash V'.
\end{align*}
This PDE can be rewritten as
\begin{align*}
\partial_t  u=\rho \bl \gamma  \sigma_\eta (\rho \nabla^2 u +2\nabla \rho \nabla u)-\kappa  W'(u)\br,\quad x\in \Omega\backslash V'.
\end{align*}
An anisotropic rescaling in time, assuming that the density $\rho\in L^\infty(\Omega)$ is bounded from below by a positive constant, yields
\begin{align*}
\partial_t  u=  \gamma  \sigma_\eta (\rho \nabla^2 u +2\nabla \rho \nabla u)-\kappa  W'(u),\quad x\in \Omega\backslash V',
\end{align*}
which is equivalent to 
\begin{align}\label{eq:gradientflow}
\partial_t (\rho u)=\rho\partial_t u=\gamma  \sigma_\eta \nabla\cdot (\rho^2 \nabla u )-\kappa \rho W'(u),\quad x\in \Omega\backslash V'.
\end{align}
Note that the anisotropic diffusion operator in \eqref{eq:gradientflow} is given by
\begin{align}\label{eq:anisodiff}
\mathcal{L}(u) = \frac{1}{\rho} \nabla \cdot (\rho^2 \nabla u)
\end{align}
which agrees with the graph Laplacian \eqref{eq:anisodiff} for the choice $(p,q,r) = (1,2,0)$. These parameters are obtained in case of thes unnormalised graph Laplacian. 
We recall that the Fiedler vector, that is the eigenvector of the second-smallest eigenvalue of \eqref{eq:anisodiff}, is used for classification. 
In our case it is not clear how the additional regularization term effects the Fiedler vector - a question which we plan to investigate in the future.

\subsection{Formal discrete to continuum limit}\label{sec:limit}
Next, we analyze the limiting behavior of weak solutions to \eqref{eq:modelweight} as $n\to \infty$ and $\epsilon\to 0$, assuming that \eqref{eq:energyinitbound} is satisfied. \\
For every test function $\phi\in C^1(\Omega)$ we have
\begin{align*}
&\frac{\di }{\di t}\int_\Omega \phi \rho_n u \di x = \frac{1}{n+m} \sum_{i=1}^{n+m} \phi(X_i)  \frac{\di}{\di t} u(X_i)\\&= \frac{1}{n+m} \sum_{i=1}^{n+m} \phi(X_i) \bl \frac{\gamma}{\epsilon^2  (n+m)} \sum_{j=1}^{ n+m} \eta_\epsilon(|X_i-X_j|) (u_j -u_i)-\kappa W'(u_i) \br\\&\quad-\frac{1}{n+m} \sum_{i=n+1}^{n+m} \phi(X_i) \bl \frac{\gamma}{\epsilon^2  (n+m)} \sum_{j=1}^{ n+m} \eta_\epsilon(|X_i-X_j|) (u_j -u_i) \br,
\end{align*}   
where we used \eqref{eq:modellabels} and that $W'(u_i)=0$ for $i=n+1,\ldots,n+m$. The second term vanishes in the limit as $n\to \infty$ since the function $u_i$ is uniformly bounded with respect  to $n$, see Proposition \ref{prop:discretesol}.
The first term can be written as 
\begin{align*}
\frac{\gamma}{\epsilon^2 }\int_\Omega\int_\Omega \phi(x) \rho_n(x) \rho_n(y) \eta_\epsilon(|x-y|) (u(y)-u(x))\di y\di x-\kappa \int_\Omega \phi(x) \rho_n(x)W'(u(x))\di x.
\end{align*}
This yields
\begin{align}\label{eq:weakformeps}
\begin{split}
\frac{\di }{\di t}\int_\Omega \phi \rho u \di x&=\lim_{n\to \infty}\frac{\di }{\di t}\int_\Omega \phi \rho_n u \di x \\
&=\frac{\gamma}{\epsilon^2 }\int_\Omega\int_\Omega \phi(x) \rho(x) \rho(y) \eta_\epsilon(|x-y|) (u(y)-u(x))\di y\di x\\
&-\kappa \int_\Omega \phi(x) \rho(x)W'(u(x))\di x,
\end{split}
\end{align} 
where we used that $\rho_n\to \rho$ weakly as $n\to\infty$ and Proposition \ref{prop:discretesol}. Note that it is sufficient to study the limit $\epsilon\to 0$ of the first term in  \eqref{eq:weakformeps} only. We use the same change of variables as in Section~\ref{sec:energyapproach} and rewrite  the first term on the right-hand side of \eqref{eq:weakformeps} as (for $\epsilon>0$ sufficiently small):
\begin{align}\label{eq:limitfirstterm}
\frac{\gamma}{\epsilon^{2}} \int_\Omega \int_{ B_r(0)} \phi(x) \rho(x) \rho(x+\epsilon z) \eta(|z|) (u(x+\epsilon z) -u(x))\di z \di x.
\end{align}
For the ease of notation, we assume that $\rho\in C^1(\Omega)$. However, the arguments can be extended to $\rho\in L^\infty(\Omega)$ by introducing a regularization of $\rho$, defined by $\rho\ast \xi_\delta$ for a non-negative, radially symmetric mollifier $\xi_\delta$ with $\delta>0$, and considering the limit as $\delta\to 0$. 
For $z=(z_1,\ldots,z_d)$, it holds
\begin{align*}
u(x+\epsilon z)-u(x)=\sum_{i=1}^d \bl u(x+\epsilon \sum_{j=1}^i z_j e_j)- u(x+\epsilon \sum_{j=1}^{i-1} z_j e_j) \br,
\end{align*}
where $e_j$ denotes the $j$th orthonormal basis vector of the $d$-dimensional Euclidean space. 
Note that any $z\in B_r(0)$ satisfies $z_j\in Q_{j-1}$ for $j=1,\ldots,d$, where $Q_{j-1}$ is the symmetric interval $Q_j=[-q_j,q_j]$ with $q_j=q_j(z_1\ldots,z_j)$ 
defined as 
$q_j=\sqrt{\bigl(r^2-\sum_{i=1}^{j} z_i^2\bigr)}.$
Then, the first term of the right-hand side of \eqref{eq:weakformeps} is given by
\begin{align}\label{eq:meanfieldhelp2}
\begin{aligned}[t]
\frac{\gamma}{\epsilon^2 }\int_\Omega&\int_\Omega \phi(x) \rho(x) \rho(y) \eta_\epsilon(|x-y|) (u(y)-u(x))\di y\di x\\
&=\frac{\gamma}{\epsilon^{2}} \sum_{i=1}^d\int_\Omega \int_{Q_0}\cdots \int_{Q_{d-1}} \phi(x) \rho(x)\rho(x+\epsilon z ) \eta(|z|)\\
&\qquad \qquad \qquad \cdot \bl u(x+\epsilon \sum_{j=1}^i z_j e_j)- u(x+\epsilon \sum_{j=1}^{i-1} z_j e_j) \br \di z_d \cdots \di z_1 \di x.
\end{aligned}
\end{align}
For the sake of readability, we consider some fixed $i\in \{1,\ldots,d\}$ and consider only the $i$th summand for now. Then, the  integral over $Q_{i-1}$ with respect to $z_i$ is given by
\begin{align}\label{eq:meanfieldhelpd}
\frac{\gamma}{\epsilon^2} \int_{Q_{i-1}} \phi(x) \rho(x)  \rho(x+\epsilon z) \eta(|z|) \bl u(x+\epsilon \sum_{j=1}^i z_j e_j)- u(x+\epsilon \sum_{j=1}^{i-1} z_j e_j) \br\di z_i.
\end{align}
Replacing $\rho(x+\epsilon z)$ by $\rho(x)+(\rho(x+\epsilon z)-\rho(x))$ allows us to write  \eqref{eq:meanfieldhelpd} as the sum of two terms. 
For the first term, we have
\begin{align*}
& \frac{\gamma}{\epsilon^2} \int_{Q_{i-1}} \phi(x) \rho^2(x)  \eta(|z|) \bl u(x+\epsilon \sum_{j=1}^i z_j e_j)- u(x+\epsilon \sum_{j=1}^{i-1} z_j e_j) \br\di z_i\\
&= {\gamma} \int_{0}^{q_{i-1}} \phi(x) \rho^2(x)\eta(|z|) z_i^2 \frac{ u(x+\epsilon \sum_{j=1}^{i-1} z_j e_j+\epsilon z_i e_i)- u(x+\epsilon \sum_{j=1}^{i-1} z_j e_j) }{\epsilon^2 z_i^2}\di z_i
\\&\quad+ \gamma\int_{0}^{q_{i-1}} \phi(x) \rho^2(x)\eta(|z|)z_i^2 \frac{ u(x+\epsilon \sum_{j=1}^{i-1} z_j e_j-\epsilon z_i e_i)- u(x+\epsilon \sum_{j=1}^{i-1} z_j e_j) }{\epsilon^2 z_i^2}\di z_i.
\end{align*}
It converges to
\begin{align*}
\gamma\int_{0}^{q_{i-1}} \phi(x) \rho^2(x)\eta(|z|) z_i^2 \partial_{x_i }^2 u(x)\di z_i
\end{align*}
in the limit $\epsilon\to 0$. 
The second term is given by 
\begin{align*}
\frac{\gamma}{\epsilon^2} \int_{Q_{i-1}} \phi(x) \rho(x) \bl( \rho(x+\epsilon z)-\rho(x)\br \eta(|z|) \bl u(x+\epsilon \sum_{j=1}^i z_j e_j)- u(x+\epsilon \sum_{j=1}^{i-1} z_j e_j) \br)\di z_i.
\end{align*}
Using
\begin{align*}
\rho(x+\epsilon z)-\rho(x)=\sum_{k=1}^d \bl \rho(x+\epsilon \sum_{l=1}^k z_l e_l)- \rho(x+\epsilon \sum_{l=1}^{k-1} z_l e_l) \br,
\end{align*}
we can write it as
\begin{align}\label{eq:meanfieldhelp}
\begin{split}
& \sum_{k=1}^d \frac{\gamma}{\epsilon^2} \int_{Q_{i-1}} \phi(x) \rho(x) \frac{ \rho(x+\epsilon \sum_{l=1}^k z_l e_l)- \rho(x+\epsilon \sum_{l=1}^{k-1} z_l e_l) }{z_k} \eta(|z|)z_k z_i \\&\qquad \qquad \qquad \frac{ u(x+\epsilon \sum_{j=1}^{i-1} z_j e_j+\epsilon z_i e_i)- u(x+\epsilon \sum_{j=1}^{i-1} z_j e_j)}{z_i} \di z_i.
\end{split}
\end{align}
Due to the symmetry of $Q_{i-1}=[-q_{i-1},q_{i-1}]$, we split the integral into  two.
Note that all summands except for  the $i$th vanish in the  limit $\epsilon\to 0$ since for   $k=1,\ldots, i-1$, we have 
\begin{align*}
&\lim_{\epsilon\to 0}\bl \frac{\gamma}{\epsilon^2} \int_{0}^{q_{i-1}} \phi(x) \rho(x) \frac{ \rho(x+\epsilon \sum_{l=1}^k z_l e_l)- \rho(x+\epsilon \sum_{l=1}^{k-1} z_l e_l) }{z_k} \eta(|z|)z_k z_i \right. \\&\qquad \qquad \qquad \bl \frac{ u(x+\epsilon \sum_{j=1}^{i-1} z_j e_j+\epsilon z_i e_i)- u(x+\epsilon \sum_{j=1}^{i-1} z_j e_j)}{z_i}\right. \\& \left. \qquad \qquad \qquad \left. \quad - \frac{ u(x+\epsilon \sum_{j=1}^{i-1} z_j e_j)- u(x+\epsilon \sum_{j=1}^{i-1} z_j e_j-\epsilon z_i e_i)}{z_i}\br \di z_i\br=0,
\end{align*}
while for  $k=i+1,\ldots,n,$ we obtain 
\begin{align*}
&\lim_{\epsilon\to 0}\bl \frac{\gamma}{\epsilon^2} \int_{0}^{q_{i-1}} \phi(x) \rho(x) \frac{ \rho(x+\epsilon \sum_{l=1}^k z_l e_l)- \rho(x+\epsilon \sum_{l=1}^{k-1} z_l e_l) }{z_k} \eta(|z|)z_k z_i \right. \\&\qquad \qquad \qquad \frac{ u(x+\epsilon \sum_{j=1}^{i-1} z_j e_j+\epsilon z_i e_i)- u(x+\epsilon \sum_{j=1}^{i-1} z_j e_j)}{z_i}
\\&\qquad \qquad - \phi(x) \rho(x) \frac{ \rho(x+\epsilon \sum_{\substack{l=1\\l\neq i}}^k z_l e_l-\epsilon z_i e_i)- \rho(x+\epsilon \sum_{\substack{l=1\\l\neq i}}^{k-1} z_l e_l-   \epsilon z_i e_i) }{z_k} \eta(|z|)z_k z_i 
\\& \qquad \qquad \qquad \left. \quad  \frac{ u(x+\epsilon \sum_{j=1}^{i-1} z_j e_j)- u(x+\epsilon \sum_{j=1}^{i-1} z_j e_j-\epsilon z_i e_i)}{z_i} \di z_i\br\\&=\gamma \int_{0}^{q_{i-1}} \phi(x)\rho(x) \bl \partial_{x_k} \rho(x)\eta(|z|) z_k z_i \partial_{x_i} u(x)-\partial_{x_k} \rho(x)\eta(|z|) z_k z_i \partial_{x_i} u(x)\br\di z_i =0.
\end{align*}
The limit of the $i$th summand in 
\eqref{eq:meanfieldhelp} as $\epsilon\to 0$ is given by
\begin{align*}
&\lim_{\epsilon\to 0}\bl \frac{\gamma}{\epsilon^2} \int_{0}^{q_{i-1}} \phi(x) \rho(x) \frac{ \rho(x+\epsilon \sum_{l=1}^i z_l e_l)- \rho(x+\epsilon \sum_{l=1}^{i-1} z_l e_l) }{z_i} \eta(|z|) z_i^2 \right. 
\\&\qquad \qquad \qquad \frac{ u(x+\epsilon \sum_{j=1}^{i-1} z_j e_j+\epsilon z_i e_i)- u(x+\epsilon \sum_{j=1}^{i-1} z_j e_j)}{z_i}
\\&\qquad \qquad + \phi(x) \rho(x) \frac{ \rho(x+\epsilon \sum_{l=1}^{i-1} z_l e_l-   \epsilon z_i e_i)- \rho(x+\epsilon \sum_{l=1}^{i-1} z_l e_l) }{z_i} \eta(|z|) z_i^2 
\\& \qquad \qquad \qquad \left. \quad  \frac{ u(x+\epsilon \sum_{j=1}^{i-1} z_j e_j-\epsilon z_i e_i)- u(x+\epsilon \sum_{j=1}^{i-1} z_j e_j)}{z_i} \di z_i\br
\\&\indent=2\gamma \int_{0}^{q_{i-1}} \phi(x)\rho(x)\partial_{x_i} \rho(x)\eta(|z|)  z_i^2 \partial_{x_i} u(x)\di z_i.
\end{align*}
Combining all of the above gives
\begin{align*}
\lim_{\epsilon \to 0}&\frac{\gamma}{\epsilon^2} \int_{Q_{i-1}} \phi(x) \rho(x)  \rho(x+\epsilon z) \eta(|z|) \bl u(x+\epsilon \sum_{j=1}^i z_j e_j)- u(x+\epsilon \sum_{j=1}^{i-1} z_j e_j) \br\di z_i
\\&=
\bl\gamma \phi(x) \rho^2(x)\partial_{x_i }^2 u(x)+2\gamma \phi(x)\rho(x)\partial_{x_i} \rho(x) \partial_{x_i} u(x) \br\int_{0}^{q_{i-1}}\eta(|z|) z_i^2 \di z_i.
\end{align*}
We still need to consider the integrals over $\Omega$ and $Q_{j-1}$ for $j\in \{1,\ldots,d\}\backslash\{i\}$ with respect to $x$ and $z_j$, respectively, 
in addition to the integral over $Q_{i-1}$. We sum over all $i=1,\ldots,d$ as in \eqref{eq:meanfieldhelp2} and obtain
\begin{align*}
\lim_{\epsilon \to 0}&\frac{\gamma}{\epsilon^2 }\int_\Omega\int_\Omega \phi(x) \rho(x) \rho(y) \eta_\epsilon(|x-y|) (u(y)-u(x))\di y\di x\\
&=
\sigma_\eta \sum_{i=1}^d\int_\Omega\gamma \phi(x) \rho^2(x)\partial_{x_i }^2 u(x)+2\gamma \phi(x)\rho(x)\partial_{x_i} \rho(x) \partial_{x_i} u(x) \di x \\
&= \sigma_\eta \left(   \gamma \int_\Omega \phi(x) \rho^2(x) \nabla^2 u(x)\di x + 2\gamma \int_\Omega \phi(x) \rho(x) \nabla \rho(x) \cdot \nabla u(x) \di x \right)\\
&= \gamma \sigma_\eta  \int_\Omega  \phi(x) \nabla \cdot (\rho^2(x) \nabla u(x))\di x,
\end{align*}
where we used that $\sigma_\eta$ is defined independently of  $i\in \{1,\ldots,d\}$  in \eqref{eq:sigmaeta} and
\begin{align*}
\int_{Q_0} \cdots \int_{Q_{i-2}}\int_{0}^{q_{i-1}}  \int_{Q_{i+1}}\cdots \int_{Q_d}\eta(|z|) z_i^2\di z_d \cdots \di z_1=\frac{1}{2} \int_{\R^d} \eta(|z|)  z_i^2 \di z =\sigma_\eta.
\end{align*}
This limiting expression allows us to state the weak formulation  \eqref{eq:weakformeps} in the   limit $\epsilon\to 0$:
\begin{align*}
\frac{\di }{\di t}\int_\Omega \phi \rho u \di x&=\gamma \sigma_\eta \int_\Omega \phi(x) \nabla \cdot(\rho^2(x) \nabla u(x)) \di x-\kappa \int_\Omega \phi(x) \rho(x)W'(u(x))\di x.
\end{align*} 
The associated mean-field PDE is given by
\begin{align*}
\partial_t (\rho u)=\rho\partial_t u=\gamma  \sigma_\eta \nabla\cdot (\rho^2 \nabla u )-\kappa \rho W'(u)
\end{align*}
which is identical to \eqref{eq:gradientflow}, derived in Section \ref{sec:energyapproach}. Note that we did not consider the boundary conditions in this derivation which are the most challenging part in terms of the rigorous derivation. A rigorous discrete to continuum result for gradient flows on graphs is proven in \cite{EPSS2019}.

\subsection{Initial and boundary conditions}
\begin{figure}[htb]
	\begin{center}
		\includegraphics[height=2.5cm,trim = {0 1.6cm 0 1.0cm}]{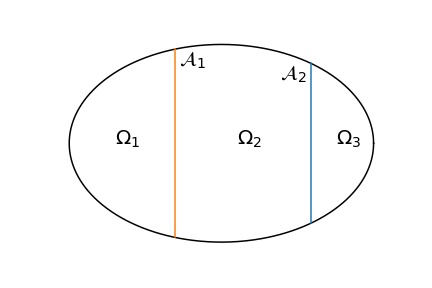}
	\end{center}
	\caption{Example of domain satisfying Assumption \ref{ass:multidim}\label{fig:domain}.}
\end{figure}

Next we discuss the corresponding initial and boundary conditions to obtain a well-posed problem.
We prescribe  initial data $u_0$ which is given by a scalar function and corresponds to the discrete interpolation of the microscopic initial condition \eqref{eq:initialdiscrete}, that is
\begin{align}\label{eq:initial}
u(x,0)=u_0(x),\quad x\in\Omega.
\end{align}
For $t>0$ and $x\in \partial \Omega\cup V'$, we extend the discrete boundary conditions \eqref{eq:bc} by considering either  Dirichlet or Neumann boundary conditions, 
as appropriate:
\begin{align}\label{eq:bccont}
u(x,t)=-1, \quad x\in V_1', \quad 
u(x,t)=1, \quad x\in V_2',\quad
\frac{\partial u}{\partial \hat{n}}(x,t)=0, \quad x\in \partial\Omega\backslash V',
\end{align}
where $\hat{n}$ denotes the exterior unit normal to the boundary $\partial \Omega$.
These boundary conditions impose restrictions on the sets $V_1'$ and $V_2'$. In 1D one can solve equation \eqref{eq:gradientflow} on intervals. The intervals are 
defined by the labeled points at which the respective boundary conditions \eqref{eq:bccont} are imposed. In higher dimensions they impose
the following additional assumptions:
\begin{ass}\label{ass:multidim}
	Suppose that there exists a $p\in\N$ and domains $\Omega_1, \ldots, \Omega_p\subset \Omega$ for $p\geq 1$ with $\Omega_i\cap \Omega_j=\emptyset$ for $i\neq j$ such 
	that $u$ is prescribed on $\mathcal A=\bar\Omega \backslash \bigcup_{i=1}^p \Omega_i\supset V'$. That is, the labels are consistent, i.e.\  there exist 
	non-empty sets $\mathcal A_1,\mathcal A_2\subset \mathcal A$ such that $\mathcal A_1\supset V'_1$ and $\mathcal A_2\supset V'_2$ with 
	$\bar{\mathcal A}_1\cap  \bar{\mathcal A}_2=\emptyset$  and $\mathcal A_1 \cup \mathcal A_2 \cup \partial \Omega = \mathcal A$ 
	such that
	\begin{align}\label{eq:bccontdimd}
	u(x,t)=-1, ~ x\in \mathcal A_1, \quad
	u(x,t)=1, ~ x\in \mathcal A_2,\quad
	\frac{\partial u}{\partial \hat{n}}(x)=0, ~ x\in \partial\Omega\backslash (\mathcal A_1\cup \mathcal A_2).
	\end{align}
\end{ass}
Note that Assumption \ref{ass:multidim} imposes necessary conditions for the continuity of $u$, in particular that
for any $x\in\mathcal A_1$ there exists an $r_1>0$ such that $B_{r_1}(x)\cap \mathcal A_2 =\emptyset$ and for any $x\in \mathcal A_2$ there  exists an 
$r_2>0$ such that $B_{r_2}(x)\cap \mathcal A_1=\emptyset$. An example of a domain, satisfying Assumption \ref{ass:multidim} is shown in Figure \ref{fig:domain}.
Assumption \ref{ass:multidim} includes the case when all labeled points are in $\partial \Omega$. Then
\eqref{eq:bccontdimd} corresponds to mixed Dirichlet and Neumann boundary conditions on $\partial \Omega$. 
For the ease of notation, we introduce the open set $\D=\bigcup_{i=1}^p \Omega_i \subset \Omega$ for $\Omega_1, \ldots, \Omega_p\subset \Omega$ in Assumption \ref{ass:multidim}.
and consider the PDE \eqref{eq:gradientflow} on $\D$ only, i.e.
\begin{align}\label{eq:gradientflowdimd}
\partial_t (\rho u)=\rho\partial_t u=\gamma  \sigma_\eta \nabla\cdot (\rho^2 \nabla u )-\kappa \rho W'(u),\quad x\in \D.
\end{align}
The respective energy \eqref{eq:energycont} is then given by
\begin{align}\label{eq:energycontdimd}
\E(u)=\frac{\gamma  \sigma_\eta}{2}\int_{\D} \rho^2 |\nabla u|^2\di x +\kappa \int_{\D}\rho W(u)\di x.
\end{align}

In summary we can formulate the continuum problem as follows: 
Let $\rho$ be a given probability measure, that is\ $\int_{\D} \rho \di x=1$, satisfying $\rho(X_i)>0$ for all $i\in V'$. 
Then continuous label propagation can described by \eqref{eq:gradientflowdimd} with initial condition \eqref{eq:initial} for a smooth 
initial label distribution obtained by the interpolation of \eqref{eq:initialdiscrete} and boundary conditions \eqref{eq:bccont}.

\subsection{Extension to multiple labels}\label{sec:multiplelabels}
We briefly comment on multi-label problems, before returning to the case of two labels in Section \ref{s:pde}.
Given $k$ labels for some $k>2$, denoted by the set of labels $L=\{l_1,\ldots,l_k\}$, we consider the set of initially labeled vertices $V'=V_1'\cup \ldots \cup V_k'$ where $V_i'\cap V_j'=\emptyset$ for $i,j\in\{1,\ldots,k\}$ with $i\neq j$. Instead of scalar-valued $u_i$, we consider $u_i\in \R^k$ where each component can take values in the interval $[0,1]$. The  label value in $[0,1]$ of the $j$th component of $u_i$ can be interpreted as the probability that $u_i$ has label $j$. The discrete boundary conditions \eqref{eq:bc} can be adapted as
\begin{align*}
u_i=e_1, \quad i\in V_1',\qquad \ldots,\qquad
u_i=e_k, \quad i\in V_k',
\end{align*}
where $e_j$ denotes the $j$th orthonormal basis vector of the $k$-dimensional Euclidean space. 
This motivates to generalize Assumption \ref{ass:multidim}  for $k$ labels:
\begin{ass}\label{ass:multidimlabels}
	Suppose that there exists a $p\in\N$ and domains $\Omega_1, \ldots, \Omega_p\subset \Omega$ for $p\geq 1$ with $\Omega_i\cap \Omega_j=\emptyset$ for $i\neq j$ such that $u\in \R^k$ is prescribed on $\mathcal A=\bar\Omega \backslash \bigcup_{i=1}^p \Omega_i\supset V'$. That is, the labels are consistent, i.e.\  there exist non-empty sets $\mathcal A_1, \ldots,\mathcal A_k\subset \mathcal A$ such that $\mathcal A_i\supset V'_i$ for $i=1,\ldots,k$ with $\bar{\mathcal A}_i\cap  \bar{\mathcal A}_j=\emptyset$ for any $i\neq j$  and $\mathcal A_1 \cup \ldots \cup \mathcal A_k \cup \partial \Omega = \mathcal A$ 
	such that
	\begin{align*}
	u(x,t)=e_i, ~ x\in \mathcal A_i, ~ 
	i=1,\ldots,k,\qquad
	J_u(x,t)\hat{n}(x,t)=0, ~ x\in \partial\Omega\backslash (\textstyle \bigcup_{i=1}^k \mathcal A_i),
	\end{align*}
	where $\hat{n}$ denotes the exterior unit normal to the boundary $\partial \Omega$ and  $J_u$ denotes the Jacobi matrix of $u$.
\end{ass}
Note that Assumption \ref{ass:multidimlabels} implies that necessary conditions for the continuity of $u$ are satisfied, i.e.\
for any $i\in\{1,\ldots,k\}$ and any $x\in\mathcal A_i$ there exists an $r_i>0$ such that $B_{r_i}(x)\cap \mathcal A_j =\emptyset$ for all $j\in\{1,\ldots,k\}\backslash \{i\}$.

\section{Properties of solutions of the continuum model}\label{s:pde}
In this section we present analytic properties of solutions to the continuum model. These results provide insights how the structure
of solutions depends on the function $\rho$ and the scaling parameters $\gamma$ and $\kappa$.\\
Throughout this section we will suppose that Assumption~\ref{ass:multidim} is satisfied and denote by 
$\D=\bigcup_{i=1}^p \Omega_i \subset \Omega$ an open set, where $\Omega_1, \ldots, \Omega_p\subset \Omega$.
Furthermore we assume: 
\begin{enumerate}
	\item[\namedlabel{itm:a1}{(A1)}] Let the density $\rho \in L^{\infty}(\D)$ with $\int_{\D} \rho \di x = 1$ and $\rho(x) \geq \bar{\rho} > 0$.
	\item[\namedlabel{itm:a2}{(A2)}] The initial datum $u_0 \in H^1(\Omega)$ satisfies $\E(u_0) < \infty$ for $\E$ given by \eqref{eq:energycontdimd}.
	\item[\namedlabel{itm:a3}{(A3)}] Let $\kappa \geq 0$.
\end{enumerate}
Further, we assume that the  double well potential $W\colon \R\to [0,\infty)$ satisfies the following assumptions:
\begin{enumerate}
	\item[\namedlabel{itm:h1}{(H1)}] $W$ is continuously differentiable,  $W(t)=0$ if and only if $t\in \{-1,1\}$.
	\item[\namedlabel{itm:h2}{(H2)}] There exist $c>0$ and $T_c>0$ such that $W(t)\geq c|t|$ for all $t\in \R$ with $|t|\geq T_c$.
\end{enumerate}
For the maximum principle, we assume in addition:
\begin{enumerate}
	\item[\namedlabel{itm:h3}{(H3)}] $W'(t)\leq 0$ for all $t\leq -1$  and $W'(t)\geq 0$ for all $t\geq 1$.
\end{enumerate}
An example for $W$ is  $W(x)=(x^2-1)^2$.

\subsection{Existence and uniqueness of solution}
\begin{theorem}\label{th:exsol}
	Let assumptions \ref{ass:multidim}, \ref{itm:a1}, \ref{itm:a2}, \ref{itm:a3}, \ref{itm:h1} and \ref{itm:h2} be satisfied. 
	Then, \eqref{eq:gradientflowdimd} subject to initial data \eqref{eq:initial} and boundary conditions \eqref{eq:bccontdimd} has a unique weak solution $u\in L^2(0,\infty;H^1(\D))$. The solution satisfies the energy dissipation inequality
	\begin{align}\label{eq:energydissipation}
	\E(u(t))+\int_0^t \int_\D (\partial_t u(x,s))^2 \di x \di s\leq \E(u_0) \quad \text{for all }t\geq 0.
	\end{align}
\end{theorem}	
For the proof of the above theorem, we introduce a regularization of \eqref{eq:gradientflowdimd}, given by the regularized PDE for $\delta>0$:
\begin{align}\label{eq:gradientflowreg}
(\rho\ast \xi_\delta) \partial_t u=\gamma  \sigma_\eta \nabla\cdot ((\rho\ast \xi_\delta)^2 \nabla u )-\kappa (\rho\ast \xi_\delta) W'(u),\quad x\in \D,
\end{align}
where $\xi_\delta$ is a non-negative, radially symmetric mollifier. Moreover, \eqref{eq:gradientflowreg} is the formal $L^2$-gradient flow of the energy
\begin{align}\label{eq:energycontreg}
\E_\delta(u)=\frac{\gamma  \sigma_\eta}{2}\int_{\D}(\rho\ast \xi_\delta)^2  |\nabla u|^2\di x +\kappa \int_{\D}(\rho\ast \xi_\delta)W(u)\di x.
\end{align} 

\begin{lemma}
	Let Assumptions \ref{ass:multidim},  \ref{itm:a1}, \ref{itm:a2}, \ref{itm:a3}, \ref{itm:h1} and \ref{itm:h2}  be satisfied and let $\delta > 0$.
	Then, the regularized PDE \eqref{eq:gradientflowreg} subject to initial data \eqref{eq:initial} and  boundary conditions \eqref{eq:bccontdimd} has a unique weak solution $u\in L^2(0,\infty;H^1(\D))$. The regularized energy \eqref{eq:energycontreg} satisfies
	\begin{align}\label{eq:energydissipationreg}
	\E_\delta(u(t))+\int_0^t \int_\D (\partial_t u(x,s))^2 \di x \di s= \E_\delta(u_0) \quad \text{for all }t\geq 0.
	\end{align}
\end{lemma}	

\begin{proof}
	The PDE \eqref{eq:gradientflowreg} with boundary conditions \eqref{eq:bccontdimd} can be regarded as multiple boundary value problems by splitting the given set $\D$ into subdomains $\Omega_1,\ldots,\Omega_p$ according to \eqref{eq:bccontdimd}.  
	Note that the sequence	$\{\rho\ast \xi_\delta\}_{\delta >0}$ is  continuously differentiable and uniformly bounded  by a positive constant from below. The term $(\rho\ast \xi_\delta) W'(u)$ is monotonically increasing in $u$ outside the compact set $[-T_c,T_c]$ by  \ref{itm:h2}.  Hence, \eqref{eq:gradientflowreg} has at least one solution global in time on $\D$ for any $\kappa\geq 0$ by the classical existence and uniqueness theory of second order parabolic PDEs \cite[§ 7.2]{evans}. The condition $\bar{\mathcal A}_1\cap  \bar{\mathcal A}_2=\emptyset$  together with the boundary conditions \eqref{eq:bccontdimd} in Assumption \ref{ass:multidim} implies the uniqueness of solutions. To obtain the energy identity \eqref{eq:energydissipationreg}, we integrate by parts
	\begin{align*}
	\frac{\di}{\di t} \E_\delta(u)&= \gamma   \sigma_\eta\int_{\D}(\rho\ast \xi_\delta)^2 \nabla u  \cdot\nabla (\partial_t u)\di x +\kappa \int_{\D}(\rho\ast \xi_\delta) \partial_t u W'(u)\di x\\&=-\int_{\D} \partial_t u  \left(\gamma   \sigma_\eta \nabla \cdot\left((\rho\ast \xi_\delta)^2 \nabla u\right) -\kappa (\rho\ast \xi_\delta)W'(u)\right)\di x=-\int_{\D} (\partial_t u)^2\di x.
	\end{align*}
	Integration in time yields \eqref{eq:energydissipationreg}.
\end{proof}

\begin{proof}[Proof of Theorem \ref{th:exsol}]
	Let $\{u_\delta\}$  be a sequence of  weak solutions of \eqref{eq:gradientflowreg}. Since $\rho$ is bounded from below by a positive constant, $\rho\ast \xi_\delta$ is also bounded from below by a positive constant, provided $\delta>0$ is chosen sufficiently small. The energy identity \eqref{eq:energydissipationreg} implies the uniform boundedness of $\|  u_\delta\|_{L^4(\D)}$ and $\| \nabla u_\delta\|_{L^2(\D)}$.
	Hence, there exists a sub-sequence, again denoted by $u_\delta$, such that we have the weak convergence $ u_\delta \rightharpoonup  u$ in $H^1(\D)$ for some $u\in H^1(\D)$ as $\delta \to 0$. 
	The strong convergence of $\rho\ast \xi_\delta$ to $\rho$ in $L^2(\D)$ implies the weak convergence  $(\rho\ast \xi_\delta) \nabla u_\delta \rightharpoonup \rho\nabla  u$. Since $H^1(\D)$ is compact in $L^2(\D)$ by the Rellich-Kondrachov theorem \cite[§ 5.7]{evans}, $u_\delta$ also converges strongly to $u$ in $L^2(\D)$.
	
	The weak formulation of \eqref{eq:gradientflowreg} is given by
	\begin{align*}
	\int_\D \phi (\rho\ast \xi_\delta) \partial_t u_\delta\di x=-\gamma \sigma_\eta \int_\D \nabla \phi \cdot ((\rho\ast \xi_\delta)^2 \nabla u_\delta )\di x-\kappa \int_\D \phi (\rho\ast \xi_\delta) W'(u_\delta)\di x
	\end{align*}
	for any test function $\phi\in H^1(\D)$. By the above convergence properties, we obtain 
	\begin{align*}
	\int_\D \phi \rho \partial_t u\di x=-\gamma \sigma_\eta \int_\D \nabla \phi \cdot (\rho^2 \nabla u )\di x-\kappa \int_\D \phi \rho W'(u)\di x
	\end{align*}
	in the limit $\delta \to 0$, i.e.\ $u$ is a weak solution of \eqref{eq:gradientflowdimd}. 
	
	The energy dissipation \eqref{eq:energydissipation} follows from passing to the limit $\delta \to 0$ in \eqref{eq:energydissipationreg} together with the weak lower semicontinuity of the $L^2$-norm and $\E_\delta(u_0)\to \E(u_0)$ as $\delta \to 0$.
\end{proof}	

\subsection{Consistency of labeling}
The PDE \eqref{eq:gradientflowdimd} is possibly strongly degenerate unless $\rho$ is bounded by a positive constant from below. 
To relax the assumptions on $\rho$, we consider the following regularization of \eqref{eq:gradientflowdimd} 
\begin{align}\label{eq:gradientflowrho}
\rho\partial_t u=\gamma  \sigma_\eta \nabla\cdot ((\varepsilon+\rho^2) \nabla u )-\kappa \rho W'(u),\quad x\in \D,
\end{align}
where $\varepsilon = \varepsilon(x) \geq \varepsilon_0 > 0$ is a prescribed function that models the  background  diffusivity. Clearly, \eqref{eq:gradientflowrho} is uniformly elliptic for any $\rho\in L^\infty(\D)$. 
\begin{corollary}
	Suppose that Assumption \ref{ass:multidim},  \ref{itm:a1}, \ref{itm:a2}, \ref{itm:a3}, \ref{itm:h1} and \ref{itm:h2}  hold.
	Then, \eqref{eq:gradientflowrho} subject to initial data \eqref{eq:initial} and boundary conditions \eqref{eq:bccontdimd} has a unique weak solution 
	$u\in L^2(0,\infty;H^1(\D))$, which satisfies the energy dissipation inequality
	\begin{align*}
	\tilde\E(u(t))+\int_0^t \int_\D (\partial_t u(x,s))^2 \di x \di s\leq \tilde\E(u_0) \quad \text{for all }t\geq 0
	\end{align*}
	with
	\begin{align}\label{eq:energycontrho}
	\tilde\E(u)=\frac{\gamma \sigma_\eta}{2}\int_{\D} (r+\rho^2) |\nabla u|^2\di x +\kappa \int_{\D}\rho W(u)\di x.
	\end{align}
\end{corollary}	
\begin{proof}
	The statement immediately follows from Theorem \ref{th:exsol}.
\end{proof}	

Note that a stationary solution $u$ to \eqref{eq:gradientflowrho} has a similar structure as the given probability measure $\rho\in L^\infty(\D)$, i.e.\ $\rho=0$ almost everywhere on some subset $\D_0\subset \D$ implies $u=0$ on $\D_0$, while $\rho>0$ and $u$ constant almost everywhere on some set $\D_c \subset \D$ implies $|u|>0$ almost everywhere on $\D_c$. 
Indeed, let $\rho\in L^\infty(\D)$ be a probability measure with $\D_0=\{x\in\D \colon \rho(x)=0 \}\neq \emptyset$. The weak formulation of \eqref{eq:gradientflowrho} implies  
\begin{align*}
\int_{\D_0} \nabla \phi \cdot  \nabla u \di x=0,
\end{align*}
for any test function $\phi\in H^1(\D)$. Therefore $\nabla u$ is zero almost everywhere on $\D_0$. In particular, $u$ is constant almost everywhere
on $\D_0$, independent of $t\geq 0$. In case of zero initial data  \eqref{eq:initial}, we have that $\rho=u=0$ on $\D_0$.
For the reverse statement, we consider the set $\tilde \D=\D\backslash \D_0$ and suppose that $u\in L^2(0,\infty;H^1(\D))$ is a weak stationary 
solution of \eqref{eq:gradientflowrho} satisfying $\nabla u=0$  almost everywhere on some set $\D_c \subset \tilde \D$, i.e.\ $u$ is constant almost everywhere on $\D_c$. 
Then, the weak formulation of \eqref{eq:gradientflowrho} implies that $W'(u)=0$ almost everywhere on $\D_c$. 
Since $\frac{\delta \tilde E}{\delta u}=\rho \partial_t u=0$ for $\tilde E$ in \eqref{eq:energycontrho}, the function $u$ is a minimizer of $\int_{\D_c} \rho W(u)\di x$ and 
therefore $u\in\{-1,1\}$ almost everywhere on $\D_c$. In particular, $u$ is non-zero almost everywhere on $\D_c$. 

We see that the function $u$ has similar characteristics as $\rho$ - a property which is desirable from a modeling point of view. It
also justifies the assumption that $\rho$ is bounded from below by a positive constant.

\subsection{Maximum principle}
Next, we show that solutions can attain their maximum and minimum on the parabolic boundary only.
For the bounded open set $\D\subset \R^d$ we define  $\D_T = \D \times (0,T]$ for some fixed time $T>0$ and $\Gamma=\{(x,t) \in \bar \D_T\colon t=0 \}$. 

\begin{theorem}\label{prop:maximumprinciple}
	Suppose that Assumption \ref{ass:multidim}, \ref{itm:h1}, \ref{itm:h2} and \ref{itm:h3} are satisfied. Let $\rho\in  L^\infty(\D)$ be a continuously differentiable probability measure bounded by a positive constant from below. Assume that $u\in C^{2,1}(\D_T)\cap C(\bar{\D}_T)$ satisfies \eqref{eq:gradientflowdimd}  subject to initial data \eqref{eq:initial} and boundary conditions \eqref{eq:bccontdimd}. Then,
	\begin{align}\label{eq:maxprinciple}
	\min_{\bar{\D}_T} u =	\min_{\Gamma\cup \mathcal A_1} u \qquad \text{and}\qquad
	\max_{\bar{\D}_T} u =	\max_{\Gamma\cup \mathcal A_2} u. 
	\end{align}
	In particular, $u(x,t)\in [-1,1]$ for all $(x,t)\in \bar \D_T$, provided  $u(x,0)\in [-1,1]$ for all $x\in\D$.  
\end{theorem}
\begin{proof}
	Consider the operator  $Lu= \gamma \sigma_\eta \nabla\cdot (\rho^2 \nabla u )$. Then, $\rho \partial_t-L$ is a uniformly parabolic operator by the assumptions on $\rho$. 
	The PDE \eqref{eq:gradientflowdimd} can be written as $\rho \partial_t u -Lu=-\kappa \rho W'(u)$. 
	For $(x,t)\in \mathcal B_1=\{(x,t)\in \D_T\colon u(x,t)<-1\}$, we have $\rho \partial_t u -Lu\geq 0$ and $\min_{\bar{\mathcal B}_1}u(x,t)$ is attained on $\bar{\mathcal B}_1\backslash \mathcal B_1=\{(x,t)\in \bar{\mathcal B}_1 \colon t=0\}\cup \{u=-1\}$ by the weak maximum principle \cite[§7.1.4, Theorem 8]{evans}. Similarly, for $\mathcal B_2=\{(x,t)\in \D_T\colon u(x,t)>1\}$, we have $\rho \partial_t u -Lu\leq 0$ and $\max_{\bar{\mathcal B}_2}u(x,t)$ is attained on $\bar{\mathcal B}_2\backslash \mathcal B_2=\{(x,t)\in \bar{\mathcal B}_2 \colon t=0\}\cup \{u=1\}$. Since $u=-1$ on $\mathcal A_1$ and $u=1$ on $\mathcal A_2$ by \eqref{eq:bccontdimd} for non-empty sets $\mathcal A_1,\mathcal A_2$, defined in Assumption \ref{ass:multidim}, we can conclude that the minimum of $u$ is attained on $\Gamma\cup \mathcal A_1$ and the maximum of $u$ is attained on $\Gamma \cup \mathcal A_2$ which yields \eqref{eq:maxprinciple}. 	
\end{proof}

Note that the requirement that $\rho$ is continuously differentiable is not restrictive for the maximum principle in Proposition \ref{prop:maximumprinciple} since properties of strong solution $u\in C^{2,1}(\D_T)\cap C(\bar{U}_T)$ of \eqref{eq:gradientflowdimd} are described.

\subsection{Dependence on edge weights and regularization}\label{sec:dependenceparam}
We conclude by studying the behavior of solutions with respect to the edge weights $w_{ij}$ and  the regularization parameter $\kappa$.\\
The definition of $\sigma_\eta$ in \eqref{eq:sigmaeta} suggests that for $\eta\leq \tilde{\eta}$ we have $\sigma_\eta\leq \sigma_{\tilde{\eta}}$. 
In particular, for $\eta(s)=\mathbbm{1}_{[0,R]}(s)$  the constant $\sigma_\eta$ increases as $R$ increases. Since  $\sigma_\eta$ only arises as a multiplicative constant in  
\eqref{eq:gradientflowdimd}, the prefactor $\gamma \sigma_\eta$ can be considered as  rescaling in time. It is therefore sufficient to investigate the dependence of $\alpha>0$ in 
\begin{align*}
\rho \partial_t u= \nabla\cdot (\rho^2 \nabla u )-\alpha \rho W'(u),\quad x\in \D,
\end{align*}
where $\alpha$ decreases as $\sigma_\eta$ increases. 
In the limit $\alpha\to 0$ its solution is equivalent to the solution of \eqref{eq:gradientflowdimd} for $\kappa=0$ up to rescaling in time. 
For $\alpha\to \infty$, or equivalently $\kappa\to \infty$ in \eqref{eq:gradientflowdimd} we obtain a $\Gamma$-convergence result (again up to a rescaling in time). In doing so
we define 
$\mathcal F_{\kappa_n}\colon L^1(\D) \to \R\cup \{+\infty\}$ for $n\in\N$ by
\begin{align*}
\mathcal F_{\kappa_n}(u)=\begin{cases}\frac{\gamma  \sigma_\eta}{2\kappa_n }\int_{\D} \rho^2 |\nabla u|^2\di x +\kappa_n \int_{\D}\rho W(u)\di x, & u\in L^4(\D)\cap H^{1}(\D), \\
+\infty, &  \text{otherwise}.
\end{cases}
\end{align*}
Here, the first case is equivalent to $\E$ in \eqref{eq:energycontdimd}  when replacing $\kappa$  by $\kappa_n$ and rescaling the functional. Minimizers of $\mathcal F_{\kappa_n}$ 
are expected to be indicator functions taking values in $\{-1,1\}$. The rescaling by $1/\kappa_n$ in the first term of $\mathcal F_{\kappa_n}$ ensures that the energy is finite.
Note that the corresponding term in the original energy \eqref{eq:energycontdimd}, blows up for all non-constant indicator functions. 

We denote by $\bv(\D;\{-1,1\})$ the space of functions of bounded variation with values in $\{-1,1\}$. For the ease of notation, we set 
$\sigma_W= \int_{-1}^{1} \sqrt{W(s)} \di s$
and introduce the weighted total variation
\begin{align*}
\text{TV}_{\rho}(u)=\sup \left\{ \int_\D u \operatorname{div}(\rho^{\frac{3}{2}} \phi) \di x\colon \phi\in C_c^\infty(\D;\R^d) \text{ s.t. } \|\phi\|_\infty\leq 1 \right\}.
\end{align*}
We define the $\Gamma$-limit $\mathcal{F}\colon L^1(\D)\to \R\cup \{+\infty\}$ of $\mathcal F_{\kappa_n}$ by
\begin{align*}
\mathcal{F}(u)=\begin{cases} \sqrt{2\gamma  \sigma_\eta}\sigma_W\text{TV}_{\rho}(u), & u\in \bv(\D;\{-1,1\}), \\
+\infty, &  \text{otherwise}.
\end{cases}
\end{align*}
This functional is motivated by the following observation:
for $u\in W^{1,2}(\D)$ and $f(t)=\int_{-1}^t \sqrt{W(s)}\di s$ with $t\in \R$ we have that $\nabla (f\circ u)=\sqrt{W(u)} \nabla u$. Therefore
\begin{align*}
\frac{\gamma  \sigma_\eta}{2\kappa_n}\int_{\D} \rho^2 | \nabla u|^2\di x +\kappa_n \int_{\D}\rho W(u)\di x &\geq 
\sqrt{2\gamma \sigma_\eta} \int_{\D} \rho^{\frac{3}{2}} |\nabla (f\circ u)| \di x\\&=\sqrt{2\gamma \sigma_\eta}\tv_\rho(f\circ u).
\end{align*}
On the other hand $u\in \bv(\D;\{-1,1\})$ implies $\tv_\rho(f\circ u)=\sigma_W\tv_\rho(u)$ since $f\circ u=f(-1)\chi_{\{u=-1\}}+f(1)\chi_{\{u=1\}}=\sigma_W\chi_{\{u=1\}}$.

\begin{theorem}
	Let $\rho\in C^1(\D)$ be a  probability measure on $\D$ which is bounded from below by a positive constant and let $\{\kappa_n\}$ be a non-negative, diverging sequence. 
	Furthermore let \ref{itm:h1} and \ref{itm:h2} be satisfied. \\ Then, $\mathcal F_{\kappa_n}$ 
	$\Gamma$-converges to $\mathcal F$ as $n\to \infty$.
	Moreover, there exists a minimizer of $\mathcal F$ in $L^1(\D)$ and 
	$$\min \{\mathcal F(u)\colon u\in L^1(\D) \}=\lim_{n \to \infty} \inf\{\mathcal F_{\kappa_n}(u)\colon u\in L^1(\D) \}.$$
	If $\{u_{\kappa_n}\}\subset L^1(\D)$ is a precompact sequence such that
	\begin{align*}
	\lim_{n\to \infty} \mathcal F_{\kappa_n}(u_{\kappa_n}) =\lim_{n\to \infty}  \inf\{\mathcal F_{\kappa_n}(u)\colon u\in L^1(\D) \},
	\end{align*}
	then every cluster point of this sequence is a minimizer of $\mathcal F$.
\end{theorem}
\begin{proof}
	By \cite[Theorem 1.21]{braides2002}, it is sufficient to show $\Gamma$-convergence and the compactness property for $\mathcal F_{\kappa_n}$, i.e.\ that any sequence $\{u_{\kappa_n}\}$ for which $\mathcal F_{\kappa_n}$ is uniformly bounded has a convergent sub-sequence. 	Since $\rho$ is bounded from above and below by positive constants, the compactness property of $\mathcal F_{\kappa_n}$ follows immediately from the compactness property of the Ginzburg-Landau functional $\mathcal F_{\epsilon}^{GL}\colon L^1(\D)\to \R$, defined as
	\begin{align}\label{eq:energyginzburglandau}
	\mathcal F_{\epsilon}^{GL}(u)=\begin{cases} \int_{\D}\epsilon|\nabla u|^2+\frac{1}{\epsilon}W(u) \di x, & u\in W^{1,2}(\D), \\
	\infty, & \text{otherwise}\end{cases}
	\end{align}
	for $\epsilon>0$ \cite{Modica1977a,Modica1977,Modica1987,Sternberg1988}. The proof of the $\Gamma$-convergence of \eqref{eq:energyginzburglandau} to 
	\begin{align*}
	\mathcal F^{GL}(u)=\begin{cases} 2\sigma_W \tv(u,\D), & u\in \bv(\D;\{-1,1\}), \\
	\infty, & \text{otherwise},\end{cases}
	\end{align*}
	as $\epsilon\to 0$, cf.\ \cite{Modica1977a,Modica1977}, can  be adapted to include the weight $\rho^{\frac{3}{2}}$.
\end{proof}

\section{Numerics}\label{s:numerics} 

We conclude by illustrating the dynamics of solutions with computational experiments. We start by discussing the respective numerical schemes on the micro- and the 
macroscopic level, before investigating the consistency across levels, the impact of the initial conditions as well as the sensitivity with respect to the model
parameters.

\subsection{Numerical methods}

\subsubsection{The microscopic solver}
The microscopic simulations are based on an explicit in time discretisation of \eqref{eq:micro}. Given the distribution of points $X_i \in \mathbb{R}^d$, $i=1, \ldots, n+m$ 
the respective weights $w_{ij}$ are computed using
\begin{align*}
w(z) = \mathbbm{1}_{\lVert z \rVert \leq R},
\end{align*}
for a suitably chosen interaction radius $R > 0$. At first we consider the case of two labels. Let
$u_i^k$ denote the label of point $X_i$ at time $t^k = k \Delta t$, where $\Delta t$ is the discrete time step. To set the initial condition $u_i^0 := u_i(0)$ we 
assign the value $1$ or $-1$ to correctly labeled points and choose one of the following three options for the remaining points:
\begin{enumerate}
	\item\label{i:homic} $u_i^0 = 0$ for all $i=1, \ldots, n$.
	\item\label{i:uniic} $u_i^0 = \mathcal{U}([-1,1])$, where $\mathcal{U}$ is the uniform distribution on $[-1,1]$.
	\item\label{i:noric} $u_i^0 = \mathcal{N}(0,\sigma^2)$, where $\mathcal{N}$ denotes the normal distribution with mean $0$ and variance $\sigma^2$.
\end{enumerate}
Note that setting the boundary conditions to a non-zero value, introduces a bias towards a specific label. While in \ref{i:homic} no prior knowledge about the 
distribution of labels is imposed, settings \ref{i:uniic} and \ref{i:noric} specify a prior weighting of labels in one direction or another.
Labels are then updated by the explicit discretisation of \eqref{eq:micro}:
\begin{align}\label{eq:microsim}
u_i^{k+1} = u_i^k  + \frac{\gamma \Delta t}{n+m} \sum_{j=1}^{n+m} w_{ij} (u_j^k - u_i^k) - \kappa \Delta t W'(u_i^k).
\end{align}
Recall that $W$ is a double well potential, which enforces binary labels $\pm 1$. Note that in contrast to the macroscopic setting
already labeled points $X_i$, $i = n, \ldots, n+m,$ can be located in the interior as well as the exterior of the computational domain.\\
In the case of $k$ labels we consider a vector valued function $u_i(t) = (u_{i,1}(t), \ldots u_{i,k}(t))$ 
instead. The components of the vector $u_i$
take values in $[0,1]$ - these values can be interpreted as the probability that the point is assigned to the respective label. Labels are updated as in \eqref{eq:microsim},
using a potential $W$ with minima at $0$ and $1$. As a natural next step one would consider an additional regularization term to enforce sparse labeling. We leave this point for future investigation. 
\subsubsection{The macroscopic solver}
The simulation of the macroscopic model \eqref{eq:gradientflowdimd} for a given point cloud requires several pre-processing steps. First the calculation of the computational
domain, which corresponds
to the convex hull of $V$. We assume that for each label at least two neighboring vertices are correctly labeled. We set Dirichlet boundary conditions on the boundary edge 
connecting the two, while assuming Neumann boundary conditions on the rest of the domain. In case of binary labels we either set $\pm 1$, while we assign a value to each label
for $k$ labels and adjust the potential $W$ accordingly. Next we calculate the respective particle density using Gaussian kernels 
\begin{align*}
\eta_\epsilon(x) = \left((2\pi)^k \lvert \Sigma \rvert\right)^{-\frac{1}{2}} e^{-\frac{1}{2} x^T\Sigma^{-1}x} \text{ with } \Sigma = \textrm{diag}(\epsilon, \ldots \epsilon) 
\end{align*}
Equation \eqref{eq:gradientflowdimd} is then solved using a splitting scheme, in which we first perform  an implicit in time step for the diffusive part
followed by an explicit in time step of the nonlinear reaction term. We use a finite difference discretization in 1D and a finite element method in 2D for each splitting step.

\subsection{Computational experiments}
For the computational experiments of the one-dimensional microscopic problem \eqref{eq:micro}  with discretization \eqref{eq:microsim},
we consider $n=250$ data points in the one-dimensional setting and parameters $\gamma=250$ and $\kappa=0.25$, unless stated otherwise.  
We generate the underlying point cloud by sampling $n/2$ points from $\mathcal{N}(-0.25, 0.125)$
and $n/2$ points from $\mathcal{N}(0.4, 0.1)$. 
For the weights, we consider $w_{ij}=\mathbbm{1}_{[0,R]}(|X_i-X_j|)$ with $R=0.25$.

\subsubsection*{Dependence on the initial data}
\begin{figure}[htbp]
	\centering
	\subfloat[Initial and stationary label distribution in case \ref{i:homic} ]{\includegraphics[scale=0.3]{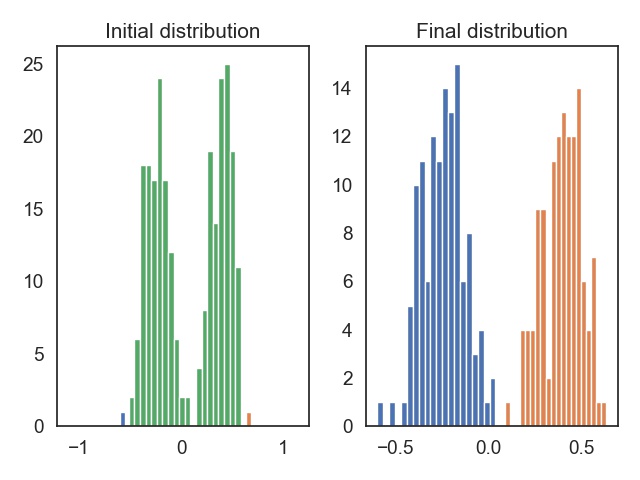}\label{fig:zeroinit}}
	\hspace*{2em}
	\subfloat[Initial and stationary label distribution in case \ref{i:noric}]{
		\includegraphics[scale=0.3]{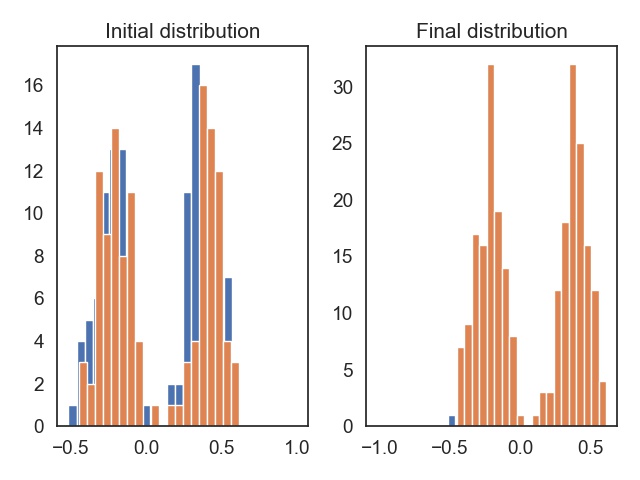}\label{fig:normone}}\\
	\subfloat[Initial and stationary label distribution in case \ref{i:uniic}]{
		\includegraphics[scale=0.3]{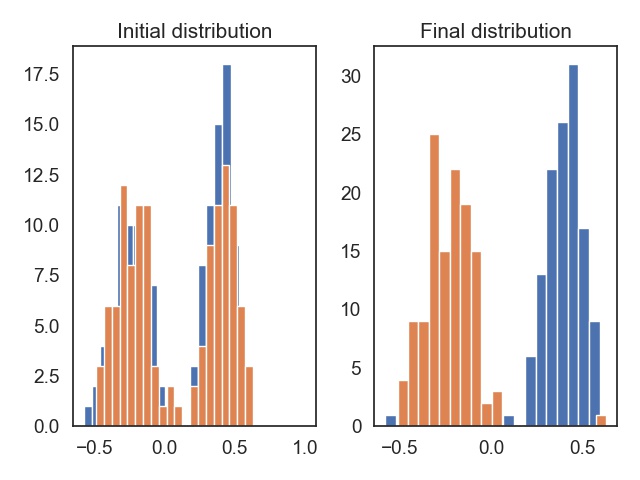}\label{fig:unimulti}}
	\hspace*{2em}
	\subfloat[Initial and stationary label distribution in case  \ref{i:uniic}]{
		\includegraphics[scale=0.3]{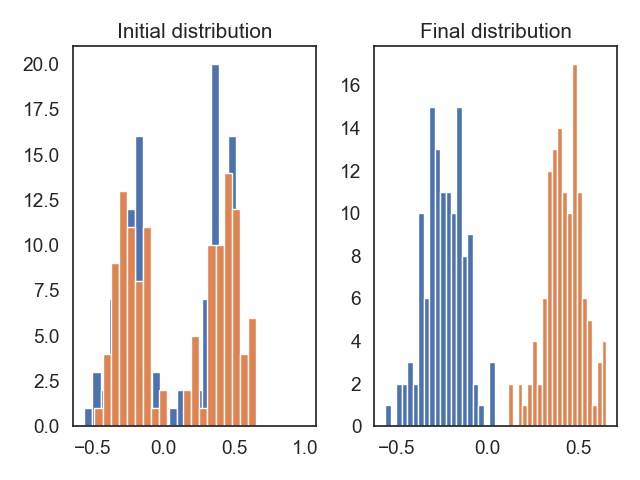}\label{fig:unitwo}}
	\caption{Dependence on initial data for different initial distributions and their stationary solution for the microscopic discretisation \eqref{eq:microsim}\label{fig:initialdata}.}
\end{figure}

To investigate the dependence on the initial data, we consider a homogeneous distribution as in case \ref{i:homic}, 
uniformly distributed points as in \ref{i:uniic} or normally distributed ones with standard deviation $0.1$ as in case~\ref{i:noric}.
The respective final label distribution is shown in Figure~\ref{fig:initialdata}, in which we use the following color coding: green corresponds to unbiased labels, that is $u_i=0$, while blue and orange encode positive or negative values of $u_i$. 
As also shown in Proposition~\ref{prop:discretesol}, the energy decays independent of the choice of initial data. While a separation into two clusters always occurs for homogeneous initial data as in Figure \subref*{fig:zeroinit}, the bias in the randomly normally or uniformly distributed case may result in  one large cluster and a  second cluster consisting only of the labeled data point at the boundary, see Figure~\subref*{fig:normone}. Moreover, the bias in the initial data may result in multiple clusters with small clusters  of the labeled points on the boundary as in Figure \subref*{fig:unimulti}. 
We have seen that the random allocation of initial labels may lead to incorrect classifications. However, if the random initial label distribution is close to the truth, then the algorithm identifies the two clusters correctly, see Figure \subref*{fig:unitwo}. This observation leads to interesting questions about the dependence of solutions on the initial datum.

\subsubsection*{Dependence of the label location}
\begin{figure}[htbp]
	\centering
	\subfloat[ Test case 1]{\includegraphics[scale=0.325]{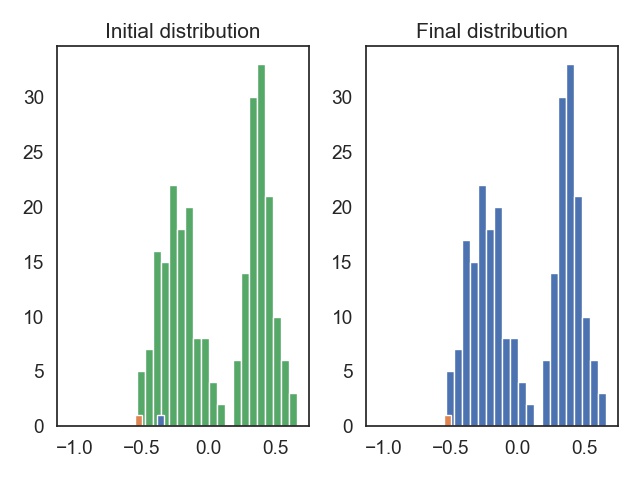}\label{fig:lablocex1}}
	\hspace*{2em}
	\subfloat[Test case 2]{
		\includegraphics[scale=0.325]{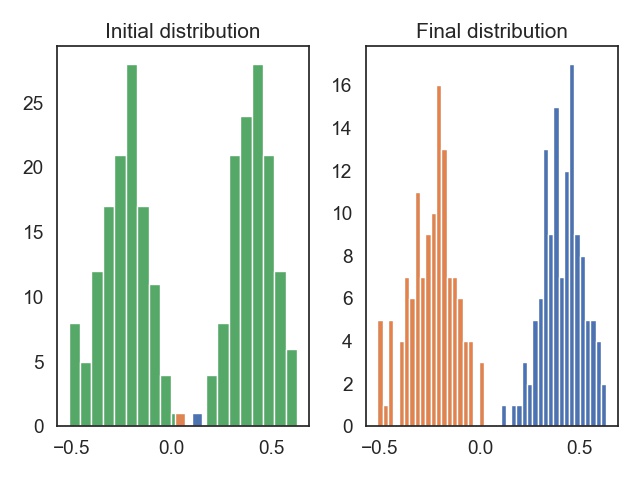}\label{fig:lablocex3}}
	\caption{Dependence on label location for initial distribution \ref{i:homic} and their stationary solution for the microscopic discretisation \eqref{eq:microsim}\label{fig:labelloc}.}
\end{figure}

For the homogeneous initial data \ref{i:homic}, we consider different label locations of the labeled data points in Figure \ref{fig:labelloc}. 
We observe that the labeling fails in the case of mislabeled points, see Figure~\subref*{fig:lablocex1}. However, if the given labels
are assigned correctly, clusters may be identified correctly even if the points are located in the interior of the point cloud, see Figure \subref*{fig:lablocex3}.  

\subsubsection*{Dependence on the scaling parameters $\gamma$ and $\kappa$}
\begin{figure}[htbp]%
	\centering
	\subfloat[Test case 1]{\includegraphics[width=0.4\textwidth]{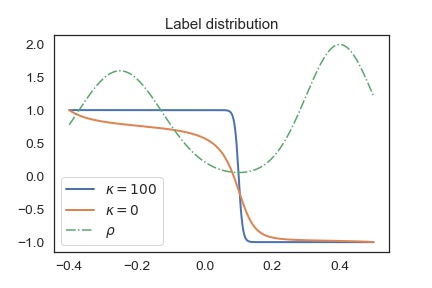}\label{fig:kappastandard}}
	\subfloat[Test case 2]{\includegraphics[width=0.4\textwidth]{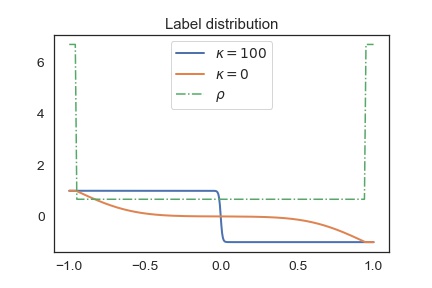}}
	\caption{Stationary solution for the macroscopic discretisation for $\kappa=0$ and $\kappa=100$ for different densities $\rho$ and homogeneous initial data \ref{i:homic}\label{fig:kappadep}.}
\end{figure}

As outlined in Section~\ref{sec:dependenceparam}, the scaling parameters $\gamma$ and $\kappa$ can be regarded as rescaling in time.
Hence it is sufficient to consider $\gamma=1$, $\sigma_\eta=\frac{1}{2}\int_{\R} \eta(|x| ) |x\cdot e|^{2} \di x=0.03125$ by \eqref{eq:sigmaeta} for $\eta(s)=\mathbbm{1}_{[0,R]}(s)$ 
with $R=0.25$  and vary $\kappa$ only. Since Section \ref{sec:dependenceparam} investigates the large data limit $n\rightarrow \infty$, we consider solutions of the macroscopic model 
\eqref{eq:gradientflowdimd} in the following. Figure \ref{fig:kappadep} illustrates the solution to \eqref{eq:gradientflowdimd} for $\kappa=0$ and $\kappa=100$ 
for given data distributions $\rho$ and homogeneous initial data. As shown in Section \ref{sec:dependenceparam}, the solution approximates an indicator function for strictly positive
functions $\rho$. Note that the solutions of the microscopic and macroscopic model are consistent due to the derivation of the limit in Section \ref{sec:limit} and numerical examples are also included in Figure \ref{fig:consistent1d}.

\subsubsection*{Consistency of the microscopic and the macroscopic model}
In the following we will investigate the consistency of the micro- and macroscopic model with 1D and 2D examples.\\

\begin{figure}[htbp]
	\centering
	\includegraphics[width=0.5\textwidth]{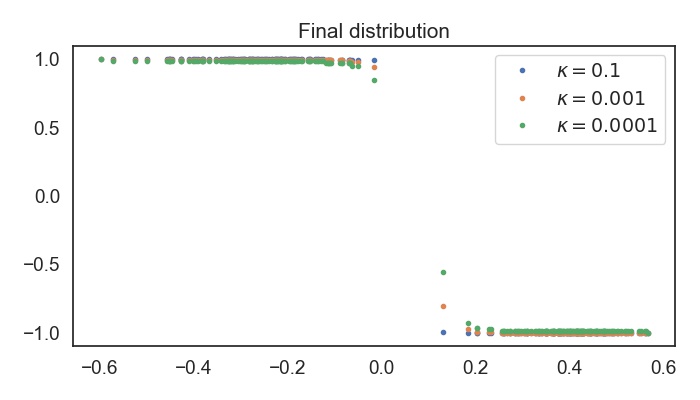}
	\caption{Stationary solution for the microscopic discretisation \eqref{eq:microsim} for different $\kappa$ and homogeneous initial data \ref{i:homic}\label{fig:consistent1d}.}
\end{figure}

\paragraph{1D}
We wish to compare the behavior of the microscopic and macroscopic solver for varying $\kappa$. In doing so we consider the original
microscopic label distribution of $250$ points equally distributed between two Gaussians. The respective density is also shown 
in Figure~\subref*{fig:kappastandard}. To compare the results of the macroscopic discretisation in Figure~\subref*{fig:kappastandard} 
to the microscopic discretisation \eqref{eq:microsim}, we consider homogeneous initial data \ref{i:homic}, set $\gamma=1$ and vary $\kappa$ in Figure~\ref{fig:consistent1d}.
We observe a very good agreement of the two levels.

\begin{figure}
	\centering
	\subfloat[Initial distribution]{\includegraphics[width=0.4\textwidth]{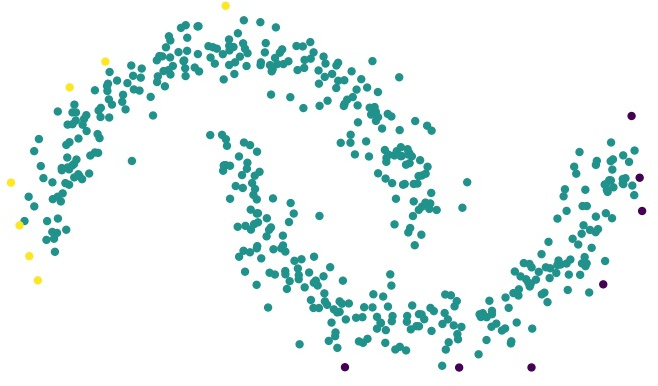}}\hspace*{3em}
	\subfloat[Final distribution]{\includegraphics[width=0.45\textwidth]{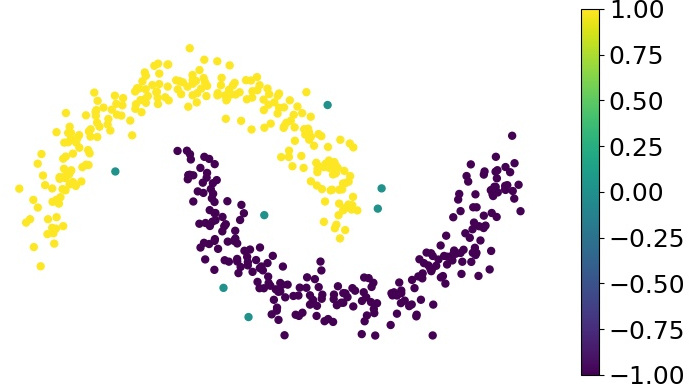}}
	\caption{Initial and final distribution at $T=25$ using the parameters $\gamma = 1.0$ and $\kappa=10$. We assume that all points defining the convex hull are correctly labelled.}
	\label{fig:twomoonsmic}
\end{figure}

\begin{figure}
	\centering
	\subfloat[Microscopic dataset]{\includegraphics[width=0.4\textwidth]{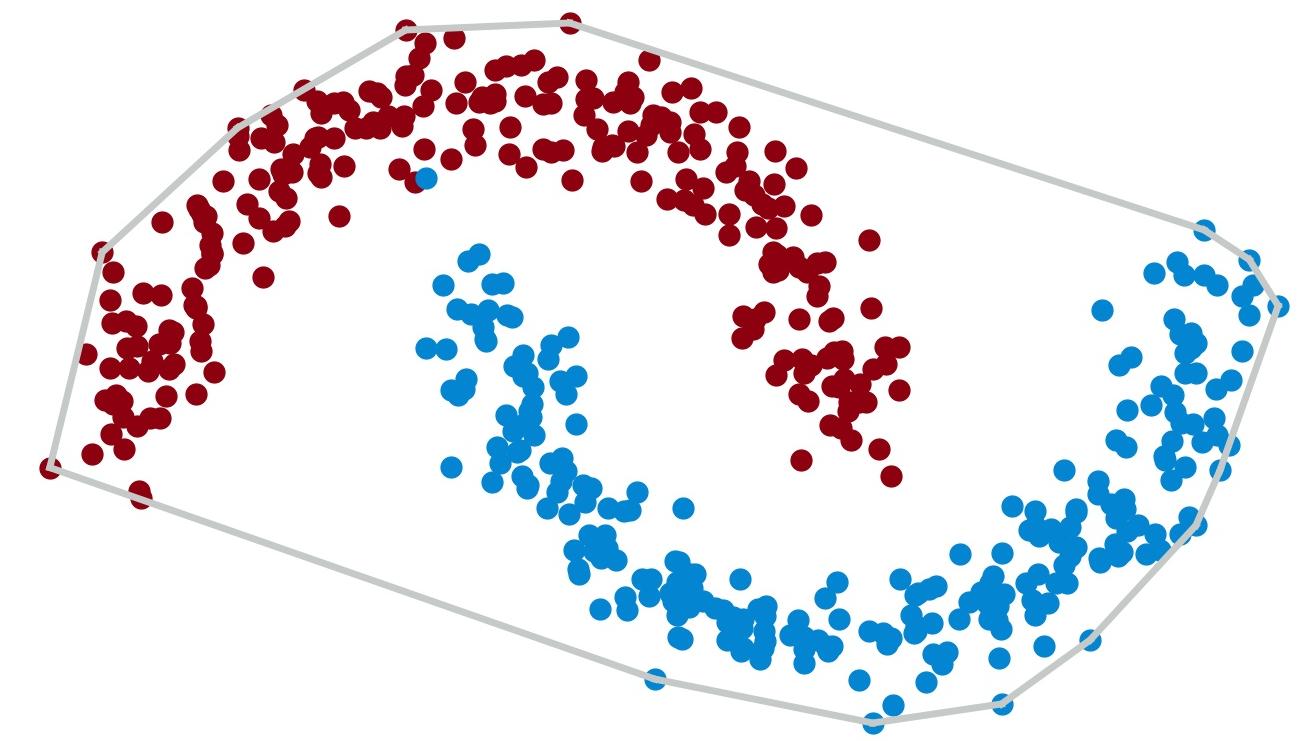}}\hspace*{3em}
	\subfloat[Final macroscopic distribution]{\includegraphics[width=0.45\textwidth]{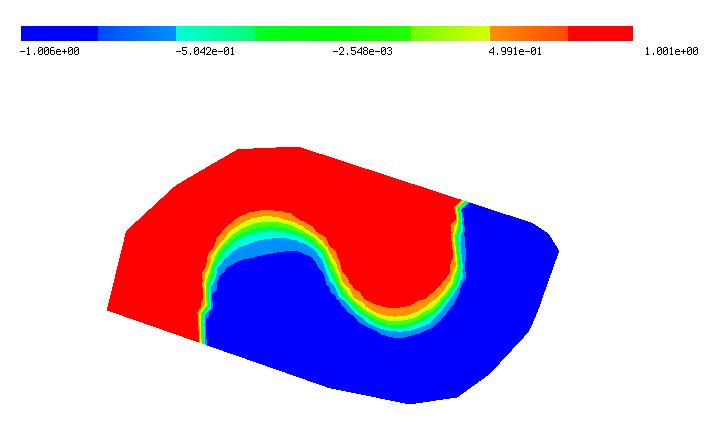}}
	\caption{Final distribution at $T=25$ using the parameters $\gamma = 1.0$ and $\kappa=10$ for the microscopic and macroscopic models. 
		We assume that all points defining the convex hull are correctly labelled.}
	\label{fig:twomoonscomp}
\end{figure}

\paragraph{Two moons} Next we investigate the performance of the microscopic and macroscopic solver for the well-known two moons problem. In doing so we generate a
noisy two moon dataset of $500$ points using Scikit, cf. \cite{scikit} with standard deviation $0.1$. We assume that the labels of points defining the convex hull
of the dataset are known and set the labels of all other points to $0$. Figure~\ref{fig:twomoonsmic} shows the initial as well as final label distribution at time $T=25$
using the parameters $\gamma = 1.0$ and $\kappa = 10$. The final distribution for the microscopic and macroscopic models in shown in Figure~\ref{fig:twomoonscomp}. 

\subsubsection*{Multilabel and more complex labeling problems} We conclude by comparing the performance of the proposed label propagation to classify single 
digits in the National Institute of Standards and Technology (NIST) dataset. In doing so we use the digits dataset of the Scikit-learn project, see 
\cite{scikit}. It consists of $1797$ samples of digital digits, which were extracted and pre-processed from the NIST dataset. Each sample corresponds to 
an $8\times 8$ image of a digit taking values in $I = \lbrace 0, 1, \ldots 16\rbrace$. We wish to assign each image to one of the numbers in $L = \lbrace 0, \ldots 9\rbrace$.\\ 
We use the  Wasserstein distance $d_{\mathcal{W}_2}$ to compare images $X_i \in I^{8\times 8}$ using the Python optimal transport (POT) library, see \cite{FC2017} and \cite{PC2020} for more 
information on computational optimal transport. 
The weights are then set to 
\begin{align*}
w_{ij} = \mathbbm{1}_{d_{\mathcal{W}_2(X_i, X_j)} \leq \bar{c}} d_{\mathcal{W}_2(X_i, X_j)}^{-1}.
\end{align*}
Images are only connected if their Wasserstein distance is below a certain threshold $\bar{c}$. Then the respective weight corresponds to the inverse of the 
Wasserstein distance. We randomly choose $320$ samples, calculate the respective weight matrix and assume that the first $40$ digits are 
correctly labelled. The remaining labels are then computed using \eqref{eq:modelorig}. To assign the images to the different labels, we use one-hot encoding - hence $u_i$ 
is a vector valued function in $\mathbb{R}^{10}$. Labels take values in the
interval $[0,1]$, which can be interpreted as the probability that the image corresponds to the respective number in $L$. 
In this case the double well potential $W$ has to be changed accordingly - with two wells located
at $x=0$ and $x=1$, in particular $W(x) = x^2(x-1)^2 $. Note that this procedure is in analogy to the extension to multiple labels in Section \ref{sec:multiplelabels}.\\
The performance of the algorithm depends strongly on the cut-off $\bar{c}$. We set it to one-tenth the value of the maximum Wasserstein distance:
\begin{align*}
\bar{c} = 0.1\times \max_{i, j \in \lbrace 1, \ldots , 320\rbrace} w_{ij}
\end{align*}
choose $\gamma=1$ and $\kappa=10$. We compare the maximum entry of each vector $u_i$, $i = 1, \ldots, 320$ at time $T=20$ with the respective correct labels - 
managing to assign $84.285\%$ of the labels correctly. We believe that the performance of the scheme can be improved by fine-tuning the parameters and adding regularization
terms to enforce sparsity of the vector $u_i$. This is however
beyond the scope of this work, in which we wish to demonstrate the general feasibility of the proposed approach.

\section{Conclusions}
In this paper we propose a continuum model for label propagation in the case of semi-supervised learning, which is formally derived from an agent based formulation.
This microscopic interpretation allows for different generalisations, for example using ideas from collective dynamics. We analyse the respective parabolic boundary value
problem for the distribution of labels, in particular the structure of these solutions. Hereby we prove that the label distribution inherits features from the 
underlying point cloud. Our computational experiments also illustrate and exemplify the connection between the structure of the underlying data set and the label propagation-
We compare the stationary distributions on the microscopic and macroscopic level.\\
The connection to agent based models allows for several generalisations - using for example powers of the absolute value $\lvert u_i - u_j \rvert^p$ or a functional
dependence $f(u_i-u_j)$. Another interesting perspective is the connection to Gaussian random fields, and Stein variational gradient descent, which were recently analysed
by  Duncan and co-workers in \cite{DNS2019}.

\section*{Acknowledgments}  LMK acknowledges support from the European Union Horizon 2020 research and innovation programmes under the Marie Sk\l odowska-Curie grant 
agreement No.\ 777826 (NoMADS) and No.\ 691070 (CHiPS), the Cantab Capital Institute for the Mathematics of Information and Magdalene College, Cambridge 
(Nevile Research Fellowship). MTW acknowledges partial support from the Austrian Academy of Sciences under the New Frontier's grant NST-0001. Both authors acknowledge support 
of the Warwick Research Development Fund through the project `Using Partial Differential Equations Techniques to Analyse Data-Rich Phenomena'.

	\bibliographystyle{plain}

	\bibliography{references}
\end{document}